\documentclass[11pt,final]{amsart}
\usepackage[utf8]{inputenc}
\usepackage[T1]{fontenc}
\usepackage[english]{babel}
\usepackage[letterpaper,centering]{geometry}
\usepackage{lmodern}
\usepackage{graphics}
\usepackage{mathrsfs}
\usepackage{amsmath,amssymb,amsfonts}

\usepackage[unicode,pdfborder={0 0 0},final]{hyperref}
\usepackage[all]{xy}
{\setbox0\hbox{$ $}}\fontdimen16\textfont2=\fontdimen17\textfont2
\entrymodifiers={+!!<0pt,\the\fontdimen22\textfont2>}
\SelectTips{cm}{11}
\usepackage{quiver}
\mathcode`A="7041 \mathcode`B="7042 \mathcode`C="7043 \mathcode`D="7044
\mathcode`E="7045 \mathcode`F="7046 \mathcode`G="7047 \mathcode`H="7048
\mathcode`I="7049 \mathcode`J="704A \mathcode`K="704B \mathcode`L="704C
\mathcode`M="704D \mathcode`N="704E \mathcode`O="704F \mathcode`P="7050
\mathcode`Q="7051 \mathcode`R="7052 \mathcode`S="7053 \mathcode`T="7054
\mathcode`U="7055 \mathcode`V="7056 \mathcode`W="7057 \mathcode`X="7058
\mathcode`Y="7059 \mathcode`Z="705A

\theoremstyle{plain}
\newtheorem{thm}{Theorem}[section]
\newtheorem{conj}[thm]{Conjecture}
\newtheorem{lem}[thm]{Lemma}

\newtheorem{cor}[thm]{Corollary}
\newtheorem*{cor*}{Corollary}
\newtheorem{prop}[thm]{Proposition}

\theoremstyle{definition}

\newtheorem{defn}[thm]{Definition}
\newtheorem{rmk}[thm]{Remark}

\numberwithin{equation}{section}

\input cyracc.def
\DeclareFontFamily{U}{russian}{}
\DeclareFontShape{U}{russian}{m}{n}
        { <5><6> wncyr5
        <7><8><9> wncyr7
        <10><10.95><12><14.4><17.28><20.74><24.88> wncyr10 }{}
\DeclareSymbolFont{Russian}{U}{russian}{m}{n}
\DeclareSymbolFontAlphabet{\mathcyr}{Russian}
\makeatletter
\let\@math@cyr\mathcyr
\renewcommand{\mathcyr}[1]{\@math@cyr{\cyracc #1}}
\makeatother
\newcommand{\Sha}{{\mathcyr{Sh}}}

\newcommand{\C}{{\mathbf C}}

\newcommand{\Q}{{\mathbf Q}}
\newcommand{\R}{{\mathbf R}}
\newcommand{\Z}{{\mathbf Z}}

\newcommand{\nr}{\mathrm{nr}}
\newcommand{\Gm}{\mathbf{G}_\mathrm{m}}
\newcommand{\et}{{\text{ét}}}

\newcommand{\SL}{\mathrm{SL}}
\newcommand{\Gal}{\mathrm{Gal}}
\newcommand{\Pic}{\mathrm{Pic}}
\newcommand{\Br}{\mathrm{Br}}
\newcommand{\Brnr}{\Br_\nr}
\newcommand{\Spec}{\mathrm{Spec}}

\newcommand{\rmH}{\mathrm{H}}
\renewcommand{\phi}{\varphi}
\renewcommand{\epsilon}{\varepsilon}

\newcommand{\Hom}{{\mathrm{Hom}}}

\newcommand{\mmu}{\boldsymbol{\mu}}

\DeclareMathOperator{\colim}{\rm colim}

\date{\today}
\title[Unramified Brauer groups of homogeneous spaces]{Unramified Brauer groups of Homogeneous spaces with Finite Stabilisers and the Grunwald problem}

\author{Lucas Lagarde}
\address{Institut Galil\'ee, Universit\'e Sorbonne Paris Nord, 99~avenue Jean-Baptiste Cl\'ement, 93430 Villetaneuse, France}
\email{lagarde@math.univ-paris13.fr}

\begin{document}
\begin{abstract} We provide an algorithm for calculating the unramified Brauer group of a homogeneous space $X$ of a semi-simple simply connected group $H$ with finite geometric stabiliser over any field of characteristic 0. When $k$ is a number field, we use the obtained description of the unramified Brauer group in order to study the Brauer-Manin obstruction to weak approximation on $X$. In particular, we provide an algorithm to compute the Brauer-Manin obstruction on $X$, which guarantees effectivity of the Grunwald problem for supersolvable groups thanks to previous work of Harpaz and Wittenberg.
\end{abstract}

\maketitle
\tableofcontents
\section{Introduction}
\label{sec:intro}

Let $k$ be a number field, let $k_v$ denote the completion of $k$ at a given place $v$, and $\Gamma_k$ denote its absolute Galois group. Let $G$ be a finite group, and let $S$ be a finite set of places of $k$. The Grunwald problem asks whether every family of prescribed local Galois extensions $L^{(v)}/k_v$ for $v\in S$, with embeddings $\Gal(L^{(v)}/k_v)\hookrightarrow G$ can be approximated by a global extension $L/k$ with Galois group $G$. More precisely:\\

\noindent\textbf{Grunwald problem.} Is the diagonal restriction map
\[\mathrm{Epi}(\Gamma_k,G)\longrightarrow\prod_{v\in S} \Hom(\Gamma_{k_v},G)/(\text{conjugacy})\]
surjective? Here $\mathrm{Epi}(\Gamma_k,G)$ denotes the set of epimorphisms $\Gamma_k\to G$ and $[P_v]= [P_v']$ in $\Hom(\Gamma_{k_v},G)/(\text{conjugacy})$ if there exists $\gamma\in G$ such that $P_v=\gamma P_v'\gamma^{-1}$.\\

This problem is of consequent interest due to its studied connections with the regular inverse Galois problem \cite{debesghazi} and with approximation properties, e.g.\ in \cite{harariquelques}. Families of groups $G$ and number fields $k$ for which $(G, k, S)$ has been studied in a systematic way and for which there is an affirmative answer to the Grunwald problem for every $S$ include abelian groups of odd order over every number field (by the Grunwald–Wang theorem), solvable groups of order prime to the number of roots of unity in $k$ (Neukirch’s theorem \cite{neukirchinv}) and groups with a generic extension over $k$ (by Saltman’s theorem \cite{saltmangeneric}).
A recent result by Lucchini Arteche \cite[Corollary 6.3]{lucchiniunramifiedbrauer} on weak approximation for homogeneous spaces specifies, for every finite group $G$, a finite set $T:= T(G,k)$ of "bad places" of $k$ such that the Grunwald problem should have an affirmative answer for $(G, k, S)$ for every set $S$ that is disjoint from $T$. The positive answer to the Grunwald problem away from the set $T$ is implied by a famous conjecture of Colliot-Thélène:

\begin{conj}\label{conjectureColliot}
Let $X$ be a smooth rationally connected variety over a number field. Then the Brauer-Manin obstruction to weak approximation on $X$ is the only one.
\end{conj}
The relation between the classical form of the Grunwald problem and this conjecture is achieved by considering the following version for a finite $k$-group:\\

\noindent\textbf{Approximation property:} A $k$-group $G$ has \emph{weak approximation} in a set $S\subset\Omega_k$ of places of $k$ if the diagonal restriction map
\[\rmH^1(k,G)\longrightarrow \prod_{v\in S}\rmH^1(k_v,G)\]
is surjective. One says that $G$ has weak approximation away from a set $T$ if $G$ has weak approximation in $S$ for every finite set $S$ of places of $k$ that is disjoint from $T$.\\

The approximation property for every finite $S\subset\Omega_k$ is equivalent to weak approximation for certain homogeneous spaces. Moreover, for constant groups $G$, it is equivalent to a positive answer to the Grunwald problem for $(G, k, S)$  for every such $S$. Recent results of Harpaz–Wittenberg \cite[Théorème B]{HWzerocycles} show that the Brauer-Manin obstruction on the associated homogeneous spaces indeed controls the approximation property, at least when the group $G(\overline{k})$ of geometric points of the $k$-group $G$ is supersolvable.

In this paper we are interested in the unramified Brauer group of homogeneous spaces of connected linear groups. In this context, when $X$
is a homogeneous space of a linear group with connected or abelian stabiliser, Borovoi, Demarche and Harari have given formulas for the "algebraic subgroup" of the unramified Brauer group $\Brnr(X)$ in \cite[Theorem~5.1]{borovoidemarcheharari}. For connected stabilisers (resp.\ finite abelian stabilisers under a semi-simple simply connected group), this result is optimal because there are no "transcendental classes" in the Brauer group, i.e.\ classes which do not vanish when passing to the algebraic closure, see \cite[Theorem 5.1]{borovoidemarcheharari} (resp.\ combine the theorem of Fischer \cite{fischer} and \cite[Proposition 26]{lucchiniartechenonramifiealgebrique}). This is not true any more when the stabiliser is non-connected, e.g.\ finite, as was already proven over the complex numbers by Saltman in \cite{saltmannoether}. Determining the whole group $\Brnr(X)$ when $k$ is arbitrary, that is, not only its algebraic subgroup, is a completely different task: a first step towards this was managed by Colliot-Thélène in \cite[Théorème 2.6]{ctnonramifie} for quotients of $\SL_{n,k}$ by finite constant groups $G$, by generalising Bogomolov’s method. Namely, he identified the normalised unramified Brauer group with a subgroup of $\Sha^2_\mathrm{ab}(G,k^*)$. More recently, Lucchini Arteche has derived in \cite[Theorem 5.1]{lucchiniunramifiedbrauer} a cohomological formula for the entire unramified Brauer group, even when the stabiliser is non-constant or if the homogeneous space $X$ has no rational point, with a small restriction on the field~$k$.

\subsection{Main results} This paper is motivated by the following well-known observation: if $X$ is a smooth, proper and rationally connected variety over a field $k$ of characteristic $0$, then the quotient group $\Br(X)/\Br_0(X)$ is finite \cite[Remarks 2.4 (ii)]{wslc}. This hints on the potential existence of an algorithm to compute the unramified Brauer group (modulo the constant classes) when $X$ is an arbitrary smooth and rationally connected variety. However the existence of such an algorithm is not known, and providing a practical one is quite out of reach in the general setting. Fixing a smooth compactification $X^c$ of $X$, the computation of the algebraic part $\Br_{\mathrm{nr},1}(X)/\Br_0(X)\hookrightarrow \rmH^1(k,\Pic(\overline{X^c}))$ would naturally rely on the existence of an algorithm to compute the Galois module structure of $\Pic(\overline{X^c})$, which is a non-trivial task when $X^c$ is rationally connected (as it already depends on the effective computability of torsion in $\ell$-adic cohomology \cite[Proposition 8.3]{poonenvanluijktesta}, \cite[Théorème~0.9]{madoreorgogozo}) and an open question of arithmetic nature when $X^c$ is arbitrary \cite[Question 2.7]{SwinnertonDyerDiophantine}. In \cite{KreschTschinkelSurface} and \cite{KreschTschinkelHassettK3}, the authors still managed to study the effectivity of the Brauer-Manin obstruction on certain families of algebraic surfaces which are not geometrically rational. However for a homogeneous space $X$ of a semi-simple simply connected group, the finiteness of $\Brnr(X)/\Br_0(X)$ is immediate since $X$ is geometrically unirational and hence rationnally connected, but there is \emph{a priori} no obvious algorithm to compute the structure of this group, let alone compute the Brauer-Manin set of $X$ when $k$ is a number field. The formula for this group provided by Lucchini Arteche in \cite[Theorem 5.1]{lucchiniunramifiedbrauer} has two main limitations: on the one hand, it has to restrict to a class of fields of characteristic zero called "non essentially real fields" (this terminology was introduced by the author in \emph{loc.\ cit}). This class is quite large (it contains all global fields and all non-archimedean local fields, for which we have natural arithmetic applications), but it excludes for instance real closed fields. On the other hand, it does not provide a framework to (practically) compute this group. In this paper, we extend this formula to arbitrary fields of characteristic zero, cf. §\ref{paragraphcritere}. In particular, we provide the following more explicit description of the unramified Brauer group when the homogeneous space admits a $k$-point:

\begin{thm}[Theorem \ref{nonramificationextension}]
Let $k$ be a field of characteristic $0$, let $G$ be a finite $k$-group together with an embedding of $k$-groups $G\hookrightarrow \SL_{n,k}$ for some $n\ge1$. Let $\alpha\in\Br^0(\SL_{n,k}/G)$ be a normalised Brauer class. Then $\alpha$ is represented by a central extension of $\Gamma_k$-groups
\[
  (E)~~: ~~ 1\longrightarrow \Q/\Z(1)\longrightarrow E(\overline{k})\longrightarrow G(\overline{k})\longrightarrow 1,
\]
and it is unramified if and only if it satisfies the following two conditions:
\begin{enumerate}
    \item[(i)] \emph{(Bogomolov condition).} The extension $(E)$ splits (as an extension of abstract groups) when pulled back to any bicyclic subgroup of $G(\overline{k})$;
    \item[(ii)] \emph{(Galois condition).} For any $\sigma\in\Gamma_k\backslash\{1\}$ and for any two elements $\tau,\gamma\in G(\overline{k})$ such that $\gamma\mkern2.5mu {}^\sigma \mkern-2.5mu\tau\gamma^{-1}=\tau^{\chi(\sigma)}$, there exists a lift of abstract groups $\psi:\langle \tau\rangle\to E(\overline{k})$ and an element $e_\gamma\in E(\overline{k})$ which maps to $\gamma$ and verifies $e_\gamma\mkern1mu {}^\sigma \mkern-1mu\psi(\tau)e_\gamma^{-1}=\psi(\tau)^{\chi(\sigma)}$.
\end{enumerate}
\end{thm}
Note that considering a homogeneous space of $\SL_{n,k}$ is not restrictive: thanks to \cite[Proposition 26]{lucchiniartechenonramifiealgebrique}, the unramified Brauer group does not depend on the choice of the ambient semi-simple simply connected group. We notably show that every geometrically unramified Brauer class descends to the base field when the latter splits the stabiliser and contains sufficiently many roots of unity:
\begin{prop}[Proposition \ref{nonramifieracinesunite}]
Let $k$ be a field of characteristic $0$ and let $G$ be a finite constant $k$-group. Fix an embedding $G\hookrightarrow\SL_{n,k}$ for some $n\ge1$ and let $N$ be the order of~$G$. Then there is an isomorphism:
\[\Brnr^0(\SL_{n,k(\zeta_{N^2})}/G)\simeq B_0(G)\]
where $B_0(G)$ denotes the Bogomolov multiplier of $G$, that is, the subgroup of $\rmH^2(G,\Q/\Z)$ consisting of classes which die under the restriction map $\rmH^2(G,\Q/\Z)\to \rmH^2(A,\Q/\Z)$ for every abelian (or equivalently, bicyclic) subgroup $A\subset G$.
\end{prop}
As a consequence, we derive a general algorithm to compute the normalised unramified Brauer group:
\begin{thm}[Theorem \ref{algorithmenonramifie}]
    There exists an algorithm that, given the datum of a field $k$ of characteristic $0$ and a finite $k$-group $G$, computes the group $\Brnr^0(\SL_{n,k}/G)$ for any embedding $G\hookrightarrow\SL_{n,k}$ for some $n\ge1$ (this does not depend on the prescribed embedding).
\end{thm}

To be precise, if $L$ denotes a finite Galois extension of $k$ over which the $k$-group $G$
 becomes constant, and if $K=L(\zeta_{N^2})$ for some primitive $(N^2)$\textsuperscript{th}-root of unity ($N$ being the order of~$G$), then the algorithm only takes as input the groups
 $\Gal(K/k)$ and $G(L)$ and the action of $\Gal(K/k)$ on both $G(L)$ and $\mmu_{N^2}\subset \Q/\Z(1)$, all of which can
 indeed be described by finitely many symbols. For more details, we refer to the proof of Theorem \ref{algorithmenonramifie}.

We also stress the fact that this algorithm computes the unramified Brauer group in the strong sense; that is, not only does it compute the structure of this group as an abelian group, but it also leads to an explicit description of the unramified classes as represented by central extensions (of $\Gamma_k$-groups) of the finite $\Gamma_k$-group $G(\overline{k})$ by the Galois module $\Q/\Z(1)$, which is well-suited with a view towards evaluation at local points when $k$ is a global field. In §\ref{paragraphkunyavskii}, we explain how to use this algorithm to construct transcendental unramified classes in the case where $k=\R$, and we give some examples.
Moreover, Theorem \ref{algorithmenonramifie} allows us to precisely describe transcendental Brauer-Manin obstructions to the Grunwald problem for finite groups. Indeed, the only known examples of such an obstruction appear in \cite{demarchelucchinineftin} and \cite{Dem22}. The first one is non-explicit: it does not tell which Brauer classes produce an obstruction. The second example computes the defect of weak approximation at local points. In §\ref{paragraphDLAN}, we determine the unramified Brauer classes that are responsible for the obstruction in the example of \cite{demarchelucchinineftin} and describe their evaluation at local points. Actually, we prove the following result:
\begin{thm}[Propositions \ref{actionsemidirect} \& \ref{obstructiontranscendante}]
Let $G=N\rtimes Q$ be a semi-direct product of abelian groups and let $k$ be a number field such that $\mmu_{\mathrm{exp}(N)}\subset k$. Then: 

\begin{enumerate} \item[(i)] the natural map $\Brnr^0(\SL_{n,k}/G)\to B_0(G)$ is surjective, and

\item[(ii)] for any transcendental Brauer class $\alpha\in\Brnr^0(\SL_{n,k}/G)$ and for any place $v\in \Omega_k$ that divides the order of $G$ and such that $Q$ can be realised as a Galois group over $k_v$, there exists a local point $P_v\in (\SL_{n,k}/G)(k_v)$ such that $\mathrm{ev}_\alpha(P_v)\neq 0$.
\end{enumerate}
In particular, any transcendental unramified Brauer class produces a Brauer-Manin obstruction to weak approximation on $\SL_{n,k}/G$.
\end{thm}
We would like to indicate that this is a strictly stronger and more precise statement than \cite[Theorem 5.1, Corollary 5.2]{demarchelucchinineftin}; indeed, the authors gave examples in \emph{loc.\ cit.}\ of semi-direct products $G=N\rtimes Q$ of abelian $p$-groups such that $\mathrm{exp}(N)=p^3$ and $\mathrm{exp}(Q)=p$, and for which there exists a transcendental Brauer-Manin obstruction to weak approximation on $\SL_{n,\Q(\zeta_{p^4})}/G$. Our approach allows to descend the transcendental obstruction to the field $\Q(\zeta_{p^3})$, which could not be achieved using their methods, see Remark \ref{remarqueDLAN}.

Coming back to arbitrary finite $k$-groups, we provide an algorithm to compute the Brauer-Manin obstruction for homogeneous spaces with finite stabiliser:
\begin{thm}[Theorem \ref{algorithmebrauermanin}]
 There exists an algorithm that, given the datum of a number field $k$ and a finite $k$-group $G$, computes the Brauer-Manin set on $\SL_{n,k}/G$ for any embedding $G\hookrightarrow\SL_{n,k}$ for some $n\ge1$.
\end{thm}

To be more precise, by "computes the Brauer-Manin set", we mean with respect to the places of $k$ which can potentially play a role in the Brauer-Manin obstruction to weak approximation on $\SL_{n,k}/G$; these places are finitely many, and we can precisely describe them thanks to \cite[Theorem~6.1, Theorem 6.2]{lucchiniunramifiedbrauer}. See the proof of Theorem \ref{algorithmebrauermanin} for more detail. Combining this result with \cite[Théorème B]{HWzerocycles}, we notably obtain that the Grunwald problem is effective for (abstract) supersolvable groups:
\begin{cor}
    There exists an algorithm that, given the datum of a finite supersolvable group $G$, a number field $k$ and a finite subset of places $S\subset \Omega_k$ provides, among the set of families, indexed by the places $v\in S$, of local Galois extensions $L^{(v)}/k_v$ with embeddings $\Gal(L^{(v)}/k_v)\hookrightarrow G$, the subset of those families that arise from a global Galois extension $L/k$ such that $\Gal(L/k)\simeq G$.
\end{cor}

\subsection{Notation and conventions} All throughout this text, if $k$ is a field, we denote by $\Gamma_k$ its absolute Galois group $\Gal(\overline{k}/k)$ where $\overline{k}$ is a fixed algebraic closure. If $k$ is a number field, then $\Omega_k$ is the set of places of $k$ and, for $v \in \Omega_k$, $k_v$ denotes the completion of $k$ at $v$. If $n\ge 1$ is an integer, we denote by $\mmu_n\subset \overline{k}^*$ the subgroup of $n$\textsuperscript{th} roots of unity in $\overline{k}$ and we write $\Q/\Z(1):=\underset{n}{\colim}\hspace{0.05cm} \mmu_n$, resp.\ $\widehat{\Z}(1):=\underset{n}{\lim}\hspace{0.05cm} \mmu_n$. Given an (abstract) group $G$, we denote by $|G|$ its order and $\mathrm{exp}(G)$ its exponent and, if $G$ is abelian and $n\ge 1$ is an integer (resp.\ $p$ is a prime integer), we denote by $G[n]$ (resp.\ $G\{p\}$) its $n$-torsion (resp.\ $p$-primary torsion) subgroup. If $G$ is finite, we let ${_\mathrm{x}}G$ be the set of subgroups $H\subseteq G$ such that $H$ is abelian ($\mathrm{x}=\mathrm{ab}$), resp.\ bicyclic ($\mathrm{x}=\mathrm{bic}$), resp.\ cyclic ($\mathrm{x}=\mathrm{cyc}$). For $i\ge 0$, if $M$ is a $G$-module, we let $\Sha^i_\mathrm{x}(G,M):=\bigcap_{H\in {_\mathrm{x}G}}\ker\big[\rmH^i(G,M)\to \rmH^i(H,M)\big]$. For instance, $B_0(G):=\Sha^2_\mathrm{bic}(G,\Q/\Z)$ is the \emph{Bogomolov multiplier} of $G$.

If $G$ is a finite $k$-group, we will make the following abuse of notation: we will always identify the finite algebraic group with the abstract group of its $\overline{k}$-points $G(\overline{k})$ seen as a $\Gamma_k$-group, i.e.\ as an abstract group equipped with a continuous action of $\Gamma_k$. For each finite $k$-group $G$, we denote by $\widehat{G}$ its Cartier dual, i.e.\ the group $\Hom_{\overline{k}}(G(\overline{k}),\overline{k}^*)\simeq \Hom_{\overline{k}}(G(\overline{k}),\Q/\Z(1))$ (morphisms taken in the category of $\overline{k}$-groups) equipped with the continuous action $(\sigma\cdot f)(g):=\sigma\cdot f(\sigma^{-1}\cdot g)$. When $k$ is algebraically closed, this group is (non-canonically) isomorphic to the Pontrjagin dual $\Hom(G,\Q/\Z)$ (we will thus adopt the same notation for both when the context is clear). A $k$-group $G$ is said to be constant if $\Gamma_k$ acts trivially on $G(\overline{k})$. 

A \emph{variety} over a field $k$ is a separated scheme of finite type over $k$. If $X$ is a variety, we denote by $\overline{X}$ the base change $X\times_k \overline{k}$ to the algebraic closure. The Brauer group of $X$ is defined as $\Br(X):=\rmH^2_\et(X,\Gm)$. This group is equipped with a filtration $\Br_0(X)\subset \Br_1(X)\subset \Br(X)$ where $\Br_0(X):=\mathrm{Im}[\Br(k)\to \Br(X)]$ along the structural morphism $X\to\Spec(k)$ and $\Br_1(X):=\ker[\Br(X)\to\Br(\overline{X})]$. If $X=H/G$ is a quotient of two $k$-groups, we denote by $\Br^0(X)$ the \emph{normalised Brauer group} of $X$, i.e. the kernel of the evaluation $1^*:\Br(X)\to\Br(k)$ at the neutral $k$-point $1\in H(k)$. The unramified Brauer group $\Brnr(X)$ is a subgroup of $\Br(X)$ which coincides with the Brauer group of any smooth compactification of $X$ (which always exist in characteristic zero thanks to Nagata's compactification theorem and Hironaka's resolution of singularities, see \cite[Proposition 6.2.7 and Corollary 6.2.10]{CTS21}). Similarly, we write $\Br_{\mathrm{nr},1}(X):=\Brnr(X)\cap \Br_1(X)$ and $\Brnr^0(X):=\Brnr(X)\cap \Br^0(X)$.

\subsection{Acknowledgement}
This work was written as part of my PhD thesis, done under the supervision of Olivier Wittenberg. I would like to express my deepest gratitude to him for his constant support, for his
patience and the time he generously gave me, as well as for his numerous suggestions and improvements on this text.
\section{Preliminaries}
\subsection{Galois covers of homogeneous spaces}
Let $X$ be a homogeneous space under a semi-simple simply connected group $H$. Let us fix a geometric point $\overline{x}$ of $X$ and assume its stabiliser $G_{\overline{x}}$ to be finite. Since $\pi_1^\et(\overline{H},\overline{x})=1$ and $\overline{X}\simeq \overline{H}/G_{\overline{x}}$, we canonically have $\pi_1^\et(\overline{X},\overline{x})=G_{\overline{x}}$. The fundamental exact sequence (see \cite[Exposé V]{SGA1})
\[1\longrightarrow \pi_1^\et(\overline{X},\overline{x})\longrightarrow \pi_1^\text{ét}(X,\overline{x})\longrightarrow\Gamma_k\longrightarrow 1\]
therefore induces a canonical outer action of $\Gamma_k$ on $G_{\overline{x}}$. Note that the isomorphism class of $G_{\overline{x}}$ seen as a group equipped with an outer action of $\Gamma_k$ does not depend on the choice of $\overline{x}$.

If $X$ admits a $k$-point $x$, there is a natural étale cover $H\to X$ given by sending $h\in H$ to $xh\in X$. If we let $G$ be the stabiliser of $x$, then $H\to X$ is a $G$-torsor which corresponds to a subgroup of $\pi_1^\text{ét}(X,\overline{x})$ which is a splitting of the projection $\pi_1^\text{ét}(X,\overline{x})\to \Gamma_k$. This yields an action of $\Gamma_k$ on $G(\overline{k})$ by conjugation in $\pi_1^\text{ét}(X,\overline{x})$, which coincides with the natural action of $\Gamma_k$ on the geometric points of the $k$-group $G$. Summing up, if $X(k)\neq\varnothing,$ then the natural extension 
\[1\longrightarrow G(\overline{k})\longrightarrow \pi_1^\text{ét}(X,\overline{x})\longrightarrow\Gamma_k\longrightarrow 1\]
splits and we may write $\pi_1^\text{ét}(X,\overline{x})\simeq G(\overline{k})\rtimes\Gamma_k$ (when $G$ is a constant $k$-group, then $\pi_1^\text{ét}(X,\overline{x})\simeq G(\overline{k})\times \Gamma_k$ is a direct product).

Let us now return to the case where $X$ does not necessarily have a $k$-point. It was remarked by Lucchini Arteche that the Brauer group of $X$ can be viewed as a cohomology group of $\pi_1^\et(X,\overline{x})$:
\newcommand{\macitation}{\cite[Proposition 3.2]{lucchiniunramifiedbrauer}}
\begin{prop}[\macitation]\label{brauerpi1}
    There exists an isomorphism
    \[\Br(X)\simeq \rmH^2(\pi_1^\et(X,\overline{x}),\Q/\Z(1))\]
    where $\pi_1^\et(X,\overline{x})$ acts on $\Q/\Z(1)$ via its quotient $\Gamma_k$ in the obvious way, and which is functorial in the following sense: for any algebraic extension $L/k$ and any homogeneous space $Y$ of $H_L$ over $L$ with finite stabiliser equipped with an $H_L$-equivariant map $Y\to X_L$, one has a commutative diagram:
\[
    \begin{tikzcd} 
          \Brnr(X) \ar[d] \ar[r, hookrightarrow] & \rmH^2(\pi_1^\et(X,\overline{x}),\Q/\Z(1))\ar[d] \\
        \Brnr(Y) \ar[r, hookrightarrow] & \rmH^2(\pi_1^\et(Y,\overline{y}),\Q/\Z(1)),
    \end{tikzcd}
\]
where the vertical arrows are the obvious restriction maps and $\overline{y}$ is any geometric point of $Y$ lying above $\overline{x}$.
\end{prop}
\subsection{Sections of the fundamental group and non-abelian cohomology} Let us keep the notation of the previous paragraph and first assume that $X$ admits a $k$-point $x$. We denote by $\mathrm{Sec}(\pi_1^\et(X,\overline{x})\to \Gamma_k)$ the set of sections of the projection $\pi_1^\et(X,\overline{x})\to \Gamma_k$. If we let $G$ be the stabiliser of any $k$-point, then there is a well-known canonical bijection of pointed sets:
\[\Sigma:\mathrm{Sec}_\sim(\pi_1^\et(X,\overline{x})\to \Gamma_k)\overset{1:1}{\longleftrightarrow} \rmH^1(k,G),\]
where $\mathrm{Sec}_\sim(\pi_1^\et(X,\overline{x})\to \Gamma_k)$ is the quotient set of $\mathrm{Sec}(\pi_1^\et(X,\overline{x})\to \Gamma_k)$ modulo conjugation by an element in $G(\overline{k})$, i.e.\ two sections $s_1,s_2:\Gamma_k\to\pi_1^\et(X,\overline{x})$ are equivalent if there exists some $g\in G(\overline{k})$ such that $s_1=gs_2 g^{-1}$. For $s\in \mathrm{Sec}(\pi_1^\et(X,\overline{x})\to \Gamma_k)$ and $\sigma\in \Gamma_k$, if we write $s(\sigma)=(h_s(\sigma),\sigma)\in \pi_1^\et(X,\overline{x})\simeq G(\overline{k})\rtimes \Gamma_k$, then $\Sigma(s)$ is nothing more than the image in $\rmH^1(k,G)$ of the cocycle $\sigma\mapsto h_s(\sigma)$. This bijection is easily seen to be compatible with field extensions.

Recall that a connected $k$-group $H$ is said to be \emph{special} if for any field extension $L/k$, we have $\rmH^1(L,H)=1$. When $X$ is a homogeneous space of $H$ which is not necessarily endowed with a $k$-point, the following result due to Pál and Schlank extends the above discussion:
\nocite*
\renewcommand{\macitation}{\cite[Theorem 9.6]{ps17}}
\begin{thm}[\macitation]\label{palschlank}
Let $H$ be a special $k$-group and $X$ a homogeneous space of $H$, let $\overline{x}$ be a geometric point of $X$ with finite stabiliser. Then for any field extension $L/k$, there is a natural map of pointed sets
\[X(L)\longrightarrow\mathrm{Sec}_\sim(\pi_1^\et(X_L,\overline{x})\to \Gamma_L)\]
which induces a bijection between $X(L)/H(L)$ and $\mathrm{Sec}_\sim(\pi_1^\et(X_L,\overline{x})\to \Gamma_L)$. This bijection is functorial in the following sense: for any intermediate field extension $k\subset E\subset L$, there is a commutative diagram of pointed sets:
\[
    \begin{tikzcd} 
          X(E)/H(E) \ar[d] \ar[r, "1:1"] & \mathrm{Sec}_\sim(\pi_1^\et(X_{E},\overline{x})\to \Gamma_E)\ar[d] \ar[l] \\
        X(L)/H(L) \ar[r, "1:1"'] & \mathrm{Sec}_\sim(\pi_1^\et(X_L,\overline{x})\to \Gamma_L), \ar[l]
    \end{tikzcd}
\]
where the vertical arrows are the obvious restriction maps.
\end{thm}

\subsection{Comparison with the Brauer-Manin pairing} Assume here that $X$ is a homogeneous space endowed with a $k$-point $x$ and let $G$ be its stabiliser, so that $\pi_1^\et(X,\overline{x})\simeq G(\overline{k})\rtimes\Gamma_k$. Let $M$ be a $\Gamma_k$-module, endowed with a $\pi_1^\et(X,\overline{x})$-action via the pullback along the projection $\pi_1^\et(X,\overline{x})\to\Gamma_k$. 
\renewcommand{\macitation}{\cite[Lemma 4.3.1]{Dem22}}
\begin{lem}[\macitation]
There exists a unique pairing of (pointed) sets:
\[(-,-):\rmH^2(\pi_1^\et(X,\overline{x}),M)\times \rmH^1(k,G)\longrightarrow \rmH^2(k,M)\]
such that $(\alpha,\Sigma(s))=s^*\alpha$ for each $\alpha\in \rmH^2(\pi_1^\et(X,\overline{x}),M)$ and $s\in\mathrm{Sec}_\sim(\pi_1^\et(X,\overline{x})\to \Gamma_k)$.
\end{lem}
\renewcommand{\macitation}{\cite[Proposition 3.4]{lucchiniunramifiedbrauer}}
\begin{lem}[\macitation]\label{comparaison}
For any field extension $L/k$, there is a commutative diagram:

\[\begin{tikzcd}
	{\rmH^2(\pi_1^\et(X_L,\overline{x}),\overline{L}^*)} & \times & {\rmH^1(L,G)} & {\rmH^2(L,\overline{L}^*)} \\
	{\Br(X_L)} & \times & {X(L)} & {\Br(L)},
	\arrow["\sim"', from=1-1, to=2-1]
	\arrow[from=1-3, to=1-4]
	\arrow[from=1-4, to=2-4]
	\arrow[from=2-3, to=1-3]
	\arrow[from=2-3, to=2-4]
\end{tikzcd}\]
where the top row is the pairing coming from the previous lemma and the bottom one is the Brauer-Manin pairing.
\end{lem}

\section{General criterion for being unramified}
\subsection{The setting} We recall that a perfect field is said to be \emph{of dimension at most $1$} if its absolute Galois group has cohomological dimension at most $1$ \cite[§II.3.1, Définition]{serrecg}. The main theorem is motivated by the following general criterion for a Brauer class to be unramified, due to Wittenberg: 
\renewcommand{\macitation}{\cite[Theorem 10.5.12]{CTS21}}
\begin{prop}[\macitation]\label{WittenbergBrnr}
Let $X$ be a smooth and connected variety over a field $k$ and let $\alpha\in \Br(X)$ be a class whose order is coprime to the characteristic exponent of $k$. Then $\alpha$ is unramified over $k$ if and only if for any perfect field $L$ of dimension at most $1$ containing $k$, the evaluation map
\[\mathrm{ev}_\alpha:X(L((t)))\longrightarrow \Br(L((t)))\]
is the zero map.   
\end{prop}
Let now $X$ be a homogeneous space of a semi-simple simply connected group $H$ over a field $k$ of characteristic $0$, with finite geometric stabiliser $G_{\overline{x}}$ (for some choice of a geometric point $\overline{x}$). The choice of some $L((t))$-point $P$ of $X$ gives rise to a commutative square:
\[
    \begin{tikzcd} 
          \Spec(L((t))) \ar[d, "P"'] \ar[dr, phantom] \ar[r] & \Spec(L)\ar[d]  \\
        X \ar[r] & \Spec(k),
    \end{tikzcd}
\]
and the covariant functoriality of the étale fundamental group induces a commutative diagram
\[
    \begin{tikzcd} 
          \Gamma_{L((t))}\ar[d] \ar[dr, phantom] \ar[r] & \Gamma_L\ar[d]  \\
        \pi_1^\text{ét}(X,\overline{x}) \ar[r, "\pi"] & \Gamma_k,
    \end{tikzcd}
\]
where $\pi:\pi_1^\et(X,\overline{x})\to\Gamma_k$ is the natural projection. Now note that since $L$ has characteristic~$0$, the group $\Gamma_{L((t))}$ can be canonically identified with the semi-direct product $\widehat{\Z}(1)\rtimes\Gamma_L$ (where the action of $\Gamma_L$ on $\widehat{\Z}(1)$ is given by the cyclotomic character): indeed, if we denote by $L((t^{\frac{1}{\infty}}))$ the field of Puiseux series in one variable with coefficients in $L$, then $\Gamma_{L((t^{\frac{1}{\infty}}))}$ can be regarded as a subgroup of $\Gamma_{L((t))}$ which projects isomorphically onto $\Gamma_L$ (this indeed holds thanks to Puiseux's theorem, see e.g.\ \cite[§XIII.2, Exercice 1]{serrecorpslocaux}). This provides a canonical section $\Gamma_L\to \Gamma_{L((t))}$ and thus $\Gamma_{L((t))}\cong \Gamma_{\overline{L}((t))}\rtimes \Gamma_L$. Moreover there is a natural isomorphism $\Gamma_{\overline{L}((t))}\simeq \widehat{\Z}(1)$, hence the above claim. 

Therefore, there is a canonically induced homomorphism $\widehat{\Z}(1)\to G_{\overline{x}}$ such that we have a commutative diagram with exact rows:
\[
    \begin{tikzcd} 
        1 \ar[r] & \widehat{\Z}(1) \ar[r] \ar[d, dashed] &  \Gamma_{L((t))} \ar[d] \ar[dr, phantom] \ar[r] & \Gamma_L\ar[d] \ar[r]& 1 \\
        1 \ar[r] & G_{\overline{x}} \ar[r]& \pi_1^\text{ét}(X,\overline{x}) \ar[r, "\pi"] & \Gamma_k \ar[r]& 1.
    \end{tikzcd}
\]

Now as it was discussed in the previous section, asking whether a class $\alpha\in\Br(X)$ evaluates trivially at $P\in X(L((t)))$ amounts to asking whether the pullback map \[P^*:\Br(X)\simeq\rmH^2(\pi_1^\et(X,\overline{x}),\Q/\Z(1))\longrightarrow \rmH^2(L((t)),\Q/\Z(1))\simeq\Br(L((t)))\] sends $\alpha$ to zero. Since $L$ is perfect of dimension at most $1$, the residue map $\Br(L((t)))\to \rmH^1(L,\Q/\Z)$ coming from the Hochschild-Serre spectral sequence
\[E_2^{p,q}=\rmH^p(L,\rmH^q(\widehat{\Z}(1),\Q/\Z(1)))\Rightarrow \rmH^{p+q}(L((t)),\Q/\Z(1))\]
is an isomorphism, see \cite[Theorem 6.3.5]{gilleszamuely}, and one has the canonical identification $\rmH^1(L,\Q/\Z)\cong \Hom_\mathrm{cont}(\Gamma_L,\Q/\Z)$. Our claim is that a careful analysis of this commutative diagram for varying choices of $L((t))$-points suffices to determine the unramified classes within the Brauer group of $X$.

\subsection{The general formula}\label{paragraphcritere} Following Lucchini Arteche \cite[Definition 2.2]{lucchiniunramifiedbrauer}, we introduce the following notation: 
\begin{defn} For $\mathrm{x}\in\{\mathrm{ab, bic, cyc}\}$ and $\mathrm{y}\in\{\mathrm{cyc}, 0\}$, we denote by ${_\mathrm{x}}(\pi_1^\et(X,\overline{x}))_\mathrm{y}$ the set of closed subgroups $D\subset \pi_1^\et(X,\overline{x})$ such that:
\begin{enumerate}
    \item[•] $D\cap G_{\overline{x}}$ is abelian $(\mathrm{x}=\mathrm{ab})$, resp.\ bicyclic $(\mathrm{x}=\mathrm{bic})$, resp.\ cyclic $(\mathrm{x}=\mathrm{cyc})$;
    \item[•] $\pi(D)\subset \Gamma_k$, is procyclic $(\mathrm{y}=\mathrm{cyc})$, resp.\ trivial $(\mathrm{y}=0)$.
\end{enumerate}
We define the following subgroups of $\Br(X)\simeq\rmH^2(\pi_1^\et(X,\overline{x}),\Q/\Z(1))$:
\[\Sha^2_{\text{x,y}}(\pi_1^\et(X,\overline{x}),\Q/\Z(1)):=\bigcap_{D\in {_\mathrm{x}}(\pi_1^\et(X,\overline{x}))_\mathrm{y}}\ker\bigg[\rmH^2(\pi_1^\et(X,\overline{x}),\Q/\Z(1))\to\rmH^2(D,\Q/\Z(1))\bigg].\]
Moreover for $\mathrm{x}\in\{\mathrm{ab, bic, cyc}\}$ we denote by $\Sha^2_{\mathrm{x},\mathrm{pcyc}}(\pi_1^\et(X,\overline{x}),\Q/\Z(1))$ the subgroup of $\rmH^2(\pi_1^\et(X,\overline{x}),\Q/\Z(1))$ consisting of those classes $\alpha$ such that, for any morphism $\widehat{\Z}\to \Gamma_k$, if we denote by $E$ the fibre product $\pi_1^\et(X,\overline{x})\times_{\Gamma_k}\widehat{\Z}$, then $\alpha$ dies in $\rmH^2(D,\Q/\Z(1))$ for any subgroup $D\subset E$ such that $D\cap G_{\overline{x}}$ is of type $\mathrm{x}$.
\end{defn}
\begin{thm}\label{criterenonramifiehomogene} Let $X$ be a homogeneous space of a semi-simple simply connected group $H$ over a field $k$ of characteristic $0$ with finite geometric stabiliser. We have the following equalities:
    \[\begin{aligned} \Brnr(X) & = \Sha^2_{\mathrm{bic},\mathrm{pcyc}}(\pi_1^\et(X,\overline{x}),\Q/\Z(1)) \\ & = \Sha^2_{\mathrm{bic},0}(\pi_1^\et(X,\overline{x}),\Q/\Z(1))\cap\Sha^2_{\mathrm{cyc},\mathrm{pcyc}}(\pi_1^\et(X,\overline{x}),\Q/\Z(1)).\end{aligned}\]

\end{thm}
\begin{rmk}
This result recovers the formula found by Lucchini Arteche in \cite[Theorem 5.1]{lucchiniunramifiedbrauer}. Our improvement lies in the fact that the proof in \emph{loc.\ cit.}\ is restricted to homogeneous spaces of semi-simple simply connected groups over non-essentially real fields, that is, fields $k$ for which the pro-$2$-Sylow subgroups of $\Gamma_k$ are either trivial or infinite. In particular, it does not apply to the field of real numbers. Our proof bypasses this restriction, see e.g.\ §\ref{paragraphkunyavskii} for examples in the case where $k=\R$.
\end{rmk}

\begin{proof} First suppose that $\alpha$ is ramified. By Proposition \ref{WittenbergBrnr} there exist a perfect field extension $L/k$ of dimension at most $1$ and some $P\in X(L((t)))$ such that $P^*\alpha\in\Br(L((t)))\simeq\Hom_\mathrm{cont}(\Gamma_L,\Q/\Z)$ is non-zero. Fix a $\widetilde{\sigma}\in\Gamma_L$ such that $P^*\alpha(\widetilde{\sigma})\neq 0$. Let $\sigma$ be the image of $\widetilde{\sigma}$ in $\Gamma_k$ along the natural morphism $\Gamma_L\to\Gamma_k$. The map
$\phi_P:\Gamma_{L((t))}=\widehat{\Z}(1) \rtimes \Gamma_L \to \pi_1^\et(X,\overline{x})$ induced
by the $L((t))$-point $P$ fits into a commutative diagram of extensions:
\[
    \begin{tikzcd} 
        1 \ar[r] & \widehat{\Z}(1) \ar[r] \ar[d, dashed] &  \widehat{\Z}(1)\rtimes\Gamma_{L} \ar[d, "\phi_P"] \ar[dr, phantom] \ar[r] & \Gamma_L\ar[d, "\tilde{\sigma}\mapsto\sigma"] \ar[r]& 1 \\
        1 \ar[r] & G_{\overline{x}} \ar[r]& \pi_1^\text{ét}(X,\overline{x}) \ar[r, "\pi"] & \Gamma_k \ar[r]& 1,
    \end{tikzcd}
\]
where the leftmost vertical map is induced by the commutativity of the right hand square and lands into a cyclic subgroup of $G_{\overline{x}}$. After pulling back the upper row of the diagram to the subgroup of $\Gamma_L$ generated by $\widetilde{\sigma}$ (which is perfect of dimension at most $1$ since $\mathrm{cd}(\Gamma_L)\le 1$, see \cite[§I.3.3, Proposition~14]{serrecg}), we can assume that $\Gamma_L$ is procyclic and that its image in $\Gamma_k$ is precisely the subgroup generated by $\sigma$. Indeed, the condition $P^*\alpha(\widetilde{\sigma})\neq 0$ is preserved thanks to the following commutative diagram provided by the functoriality of the Hochschild-Serre spectral sequence:
\[
    \begin{tikzcd} 
   \Br(\overline{L}^{\tilde{\sigma}}((t))) \ar[r, "\sim"] & \Hom_\mathrm{cont}(\Gamma_{\overline{L}^{\tilde{\sigma}}},\Q/\Z)   \\
        \Br(L((t))) \ar[u] \ar[r, "\sim"] & \Hom_\mathrm{cont}(\Gamma_L,\Q/\Z).\ar[u] 
    \end{tikzcd}
\]

Let $\zeta\in \widehat{\Z}(1)$ be a topological generator and denote by $(\zeta,1)$ its image in $\widehat{\Z}(1)\rtimes\Gamma_L$. Let $D\subset\pi_1^\et(X,\overline{x})$ be the image of $\varphi_P$. We see that $D\cap G_{\overline{x}}$ is generated by $\tau:=\varphi_P(\zeta,1)\in G_{\overline{x}}$ and by the image $\gamma$ of some topological generator $\sigma'$ of $\ker[\Gamma_L\to\Gamma_k]$. The semidirect product structure on $\widehat{\Z}(1)\rtimes \Gamma_L$ is given by the formula $(1,\sigma)(\zeta,1)(1,\sigma)^{-1}=(\zeta^{\chi(\sigma)},1)$ (where $\chi:\Gamma_L\to \mathrm{Aut}(\widehat{\Z}(1))=\widehat{\Z}^*$ is the cyclotomic character). By the definition of $\varphi_P$ and since $\sigma'$ dies in $\Gamma_k$, we obtain in particular that $\gamma\tau\gamma^{-1}=\tau$, so that $D\subset \pi_1^\et(X,\overline{x})$ is a subgroup such that $D\cap G_{\overline{x}}$ is bicyclic and $\pi(D)=\langle\sigma\rangle\subset\Gamma_k$. 

Remark that $\alpha$ does not die when pulled back to $\rmH^2(D,\Q/\Z(1))$ as it does not when pulled back to $\Br(L((t)))$. Actually, since $\Gamma_L$ is procyclic of cohomological dimension at most $1$, one can extend $\Gamma_L\to \pi(D)$ to a surjection $\widehat{\Z}\twoheadrightarrow \pi(D)$ by considering $\widehat{\Z}\to \Gamma_L$ sending $1\mapsto \widetilde{\sigma}$. Since $\alpha$ survives in $\Hom(\Gamma_L,\Q/\Z)$, \emph{a fortiori} in $\Hom(\widehat{\Z},\Q/\Z)$, then it survives in $\rmH^2(E,\Q/\Z(1))$ where $E$ denotes the fibre product $D\times_{\pi(D)}\widehat{\Z}$. We thus obtain the inclusion $\Sha^2_{\mathrm{bic},\mathrm{pcyc}}(\pi_1^\et(X,\overline{x}),\Q/\Z(1))\subset \Brnr(X)$.

Let us now carefully analyse the structure of the group $D$. If $\sigma=1\in \Gamma_k$, then $D=D\cap G_{\overline{x}}$ is bicyclic without a prescribed Galois action and $\alpha$ survives in $\rmH^2(D,\Q/\Z(1))$, that is, $\alpha$ does not lie in $\Sha^2_{\mathrm{bic},0}(\pi_1^\et(X,\overline{x}),\Q/\Z(1))$. We may now assume that $\alpha\in\Sha^2_{\mathrm{bic},0}(\pi_1^\et(X,\overline{x}),\Q/\Z(1))$ and that $\sigma\neq 1$. Observe that the Galois structure on $B:=D\cap G_{\overline{x}}$ is compatible with the inclusions $\langle \tau\rangle \subset B$ and $\langle \gamma\rangle\subset B$ in the sense that both of these cyclic subgroups are stable under the action of $\mathrm{Im}[\Gamma_L\to \Gamma_k]$: indeed we have, by the very construction of $D$, a commutative diagram with exact rows: 
\[
    \begin{tikzcd} 
        1 \ar[r] & \widehat{\Z}(1)\times \ker[\Gamma_L\to\Gamma_k] \ar[r] \ar[d] &  \widehat{\Z}(1)\rtimes\Gamma_{L} \ar[d, "\varphi_P"] \ar[dr, phantom] \ar[r] & \mathrm{Im}[\Gamma_L\to\Gamma_k]\ar[d, hookrightarrow] \ar[r]& 1 \\
        1 \ar[r] & G_{\overline{x}} \ar[r]& \pi_1^\text{ét}(X,\overline{x}) \ar[r, "\pi"] & \Gamma_k \ar[r]& 1,
    \end{tikzcd}
\]
so that both $\langle \tau\rangle$ and $\langle \gamma\rangle$  are $\mathrm{Im}[\Gamma_L\to \Gamma_k]$-stable, as desired. Let $K:=\overline{k}^{\sigma}$ and let $Y$ be a homogeneous space of $H_K$ with geometric stabiliser isomorphic to $B=D\cap G_{\overline{x}}$ together with a morphism $\psi:Y\to X$ (the existence of such a morphism is implied by e.g.\ \cite[Proposition 3.1]{lucchiniunramifiedbrauer} applied to $D\subset \pi_1^\et(X,\overline{x})$). Since $\alpha\in\Sha^2_{\mathrm{bic},0}(\pi_1^\et(X,\overline{x}),\Q/\Z(1))$, we see that $\psi^*(\alpha)$ dies in $\Br(\overline{Y})\simeq\rmH^2(B,\Q/\Z)$, hence is an algebraic Brauer class, and similarly for its image in $\Br(Y_L)$. The Hochschild-Serre spectral sequence for $Y_L\to \Spec(L)$ gives rise to a short exact sequence:
\[\Br(L)\longrightarrow \ker[\Br(Y_L)\to\Br(\overline{Y_L})]\longrightarrow\rmH^1(L,\widehat{B})\longrightarrow 0,\]
and $\Br(L)=0$ (as $\Gamma_L$ has cohomological dimension $1$). Note that the fundamental exact sequence 
\[1\longrightarrow B\longrightarrow \pi_1^\et(Y_L,\overline{y})\longrightarrow \Gamma_L\longrightarrow 1\]
splits because $\Gamma_L$ has cohomological dimension $1$ (here $\overline{y}$ is any geometric point of $Y_L$ lying above $\overline{x}$). Let us therefore fix an isomorphism $\pi_1^\et(Y_L,\overline{y})\simeq B\rtimes\Gamma_L$. Since both $\langle \tau\rangle$ and $\langle\gamma\rangle$ are $\Gamma_L$-stable,
one can define morphisms of semidirect products $\langle\tau\rangle\rtimes\Gamma_L\to B\rtimes\Gamma_L$ and $\langle\gamma\rangle\rtimes\Gamma_L\to B\rtimes\Gamma_L$ in the obvious way. Applying \cite[Proposition 3.1]{lucchiniunramifiedbrauer} to $\langle\tau\rangle\rtimes\Gamma_L$ and $\langle\gamma\rangle\rtimes\Gamma_L$ viewed as subgroups of $\pi_1^\et(Y_L,\overline{y})$ defines two homogeneous spaces $W$ and $Z$ of $H_L$ with geometric stabiliser $\langle \tau\rangle$ and $\langle\gamma\rangle$ respectively, together with morphisms $W\to Y_L$ and $Z\to Y_L$. By assumption the image $\beta$ of $\psi^*(\alpha)$ in $\Br_1(Y_L)\simeq \rmH^1(L,\widehat{B})$ is non-zero, hence by functoriality of the Hochschild-Serre spectral sequence it cannot simultaneously die in both $\Br_1(W)\simeq\rmH^1(L,\widehat{\langle\tau\rangle})$ and $\Br_1(Z)\simeq\rmH^1(L,\widehat{\langle\gamma\rangle})$. This proves the inclusion
\[\Sha^2_{\mathrm{bic},0}(\pi_1^\et(X,\overline{x}),\Q/\Z(1))\cap\Sha^2_{\mathrm{cyc},\mathrm{pcyc}}(\pi_1^\et(X,\overline{x}),\Q/\Z(1))\subset \Brnr(X).\]
Now suppose $\alpha$ is unramified. Fix a procyclic subgroup $\Gamma\subset \Gamma_k$ together with a surjection $\widehat{\Z}\to \Gamma$, and consider the fibre product $E:=\pi_1^\et(X,\overline{x})\times_\Gamma \widehat{\Z}$. The natural map $E\to \pi_1^\et(X,\overline{x})$ factors through $\pi_1^\et(X_K,\overline{x})$ where $K$ is the subfield of $\overline{k}$ fixed by $\Gamma$, and since $\widehat{\Z}$ is a free profinite group we see that $E$ is isomorphic to the semidirect product $G_{\overline{x}}\rtimes\widehat{\Z}$ (where $\widehat{\Z}$ acts on $G_{\overline{x}}$ through its quotient $\Gamma$). If we fix a subgroup $D\subset E$ such that $D\cap G_{\overline{x}}$ is bicyclic, then the image of $D$ in $\pi_1^\et(X_K,\overline{x})$ is a subgroup $D'$ such that $D'\cap G_{\overline{x}}$ is abelian (since it is generated by the images of the bicyclic group $D\cap G_{\overline{x}}$ and the procyclic group $\ker[\widehat{\Z}\to \Gamma_K]$ which commute). The pullback $\Br(X)\to \rmH^2(D',\Q/\Z(1))$ identifies to a morphism $\Br(X)\to \Br(Y)$ where $Y$ is a homogeneous space of $H_K$ with geometric stabiliser isomorphic to $D'\cap G_{\overline{x}}$, and it moreover restricts to a morphism $\Brnr(X)\to\Brnr(Y)$ thanks to Proposition \ref{brauerpi1}. By a result of Lucchini Arteche \cite[Proposition~4.1]{lucchiniunramifiedbrauer}, and as $\Gamma_K$ is procyclic and $D'\cap G_{\overline{x}}$ is abelian, we have $\Brnr(Y)\simeq \Br_0(Y)$. Now note that we have a commutative square:
\[
    \begin{tikzcd} 
   \rmH^2(\widehat{\Z},\Q/\Z(1)) \ar[r] & \rmH^2(D,\Q/\Z(1))   \\
        \Br(K) \ar[u] \ar[r] & \rmH^2(D',\Q/\Z(1))\ar[u],
    \end{tikzcd}
\]
where $\rmH^2(\widehat{\Z},\Q/\Z(1))$ is zero. The image of $\alpha$ in $\rmH^2(D,\Q/\Z(1))$ must therefore be zero since it comes from a constant. Hence $\alpha$ belongs to $\Sha^2_{\mathrm{bic},\mathrm{pcyc}}(\pi_1^\et(X,\overline{x}),\Q/\Z(1))$, and the inclusion
\[\Sha^2_{\mathrm{bic},\mathrm{pcyc}}(\pi_1^\et(X,\overline{x}),\Q/\Z(1))\subset\Sha^2_{\mathrm{bic},0}(\pi_1^\et(X,\overline{x}),\Q/\Z(1))\cap\Sha^2_{\mathrm{cyc},\mathrm{pcyc}}(\pi_1^\et(X,\overline{x}),\Q/\Z(1))\]
implies the statement of the theorem.
\end{proof}
\begin{cor}\label{unramifiednormalisehomogene}
Keeping the same notation, if $X(k)\neq\varnothing$, then for any geometric point $\overline{x}$ lying over a $k$-point $x$ of $X$, we have the equality $\Br_\mathrm{nr}(X)=\Brnr^0(X)\oplus \Br(k)$, where
    \[\Brnr^0(X)=\bigg(\Sha^2_{\mathrm{bic},0}(G_{\overline{x}}\rtimes\Gamma_k,\Q/\Z(1))\cap\Sha^2_{\mathrm{cyc},\mathrm{pcyc}}(G_{\overline{x}}\rtimes\Gamma_k,\Q/\Z(1))\bigg)\cap\Br^0(X).\]

\end{cor}
\begin{proof}
If we assume $X(k)$ to be non-empty, then $X\simeq H/G$ for some finite $k$-group $G$ such that $G(\overline{k})\simeq G_{\overline{x}}$ (as $\Gamma_k$-groups). The $k$-point on $X$ given by the image of $1\in H(k)$ provides a splitting $\Br(X)=\Br^0(X)\oplus \Br(k)$ (where both $\Br^0(X)\simeq \Br(X)/\Br(k)$ and $\Br(k)$ are subgroups of $\Br(X)$), as well as an isomorphism $\pi_1^\et(X,\overline{x})\simeq \pi_1^\et(X,\overline{1})\simeq G(\overline{k})\rtimes\Gamma_k\simeq G_{\overline{x}}\rtimes \Gamma_k$. Applying Theorem \ref{criterenonramifiehomogene} provides the result.
\end{proof}
\begin{cor}\label{algorithmecorpsreelsclos}
   Let $X$ be a homogeneous space of a semi-simple simply connected group over a real closed field $k$, with finite geometric stabiliser. Then the group $\Brnr(X)$ is effectively computable.
\end{cor}
\begin{proof}
    Indeed, if $\overline{x}$ is a geometric point of $X$, then $\pi_1^\et(X,\overline{x})$ is a finite group as an extension of $\Gamma_k\simeq \Z/2$ by the stabiliser $G_{\overline{x}}$ of $\overline{x}$, so that $\Br(X)\simeq\rmH^2(\pi_1^\et(X,\overline{x}),\Q/\Z(1))$ already is a finite group, whose elements are classes of extensions of $\pi_1^\et(X,\overline{x})$ by $\Q/\Z(1)$ (endowed with its action of $\pi_1^\et(X,\overline{x})$ through the quotient $\Gamma_k$). Checking whether a Brauer class lies in the groups $\Sha^2_{\mathrm{bic},0}(\pi_1^\et(X,\overline{x}),\Q/\Z(1))$ and $\Sha^2_{\mathrm{cyc},\mathrm{pcyc}}(\pi_1^\et(X,\overline{x}),\Q/\Z(1))$ is done by looking at its restriction with respect to subgroups of $\pi_1^\et(X,\overline{x})$, so that one can conduct the verification that this class is unramified in a finite number of steps.
\end{proof}
\section{Refined criteria for being unramified}
\subsection{General identifications} We first provide a more concrete description of the Brauer group of a  homogeneous space with finite stabiliser together with a rational point. Let us start by recalling the following well-known fact:
\renewcommand{\macitation}{\cite[Proposition 2.10]{colliotthélènexu}}
   \begin{lem}[\macitation]
    Let $H$ be a semi-simple simply connected group, let $G$ be a $k$-subgroup of $H$ and $\phi:H\to H/G$ the associated $G$-torsor. Then the natural morphism $\phi^*:\Br(H/G)\to \Br(H)$ identifies to the evaluation $1^*:\Br(H)\to \Br(k)$ at the neutral point $1\in H(k)$, and its kernel is canonically isomorphic to $\Br(H/G)/\Br(k)$.
\end{lem}
\begin{defn} Let $\Gamma$ be a profinite group and $G$ be a finite group endowed with a continuous action of $\Gamma$. If $M$ is a discrete $\Gamma$-module, we denote by $\mathrm{Ext}_\Gamma^c(G,M)$ the group of (classes of) central extensions, in the category of $\Gamma$-groups, of $G$ (equipped with its action of~$\Gamma$) by $M$. In other words, we require for all the maps appearing in these central extensions to be $\Gamma$-equivariant. 
\end{defn}

When $\Gamma=\Gamma_k$ is the absolute Galois group of a field $k$ and $G$ is a finite $k$-group, we will instead denote $\mathrm{Ext}_{\Gamma_k}^c(G(\overline{k}),M)$ by $\mathrm{Ext}_k^c(G(\overline{k}),M)$ for simplicity. The following result is inspired by \cite[Proposition 1.1]{demarchenilpotent} (the given reference treats the case where the $k$-group $G$ is constant):
\begin{prop}\label{normaliseextension}
    Let $G$ be a finite $k$-group together with an embedding into a semi-simple simply connected algebraic $k$-group $H$. Then there exists an isomorphism
    \[\mathrm{Ext}^c_k(G(\overline{k}),\Q/\Z(1))\overset{\sim}{\longrightarrow}\Br^0(H/G),\]
    which is functorial in the $k$-group $G$.
\end{prop}
\begin{proof}
  By Proposition \ref{brauerpi1}, we have $\Br(H/G)\simeq \rmH^2(\pi_1^\et(H/G,\overline{1}),\Q/\Z(1))$ where $1\in (H/G)(k)$ denotes the image of the neutral point $1\in H(k)$ (viewed as a geometric point of $H/G$), and $\pi_1^\et(H/G,\overline{1})$ acts on $\Q/\Z(1)$ through its quotient $\Gamma_k$. For any integer $n\ge 1$, the Kummer sequence on $H/G$ induces a short exact sequence
\[0\longrightarrow \Pic(H/G)/n\longrightarrow \rmH^2(\pi_1^\et(H/G,\overline{1}),\mmu_n)\longrightarrow \Br(H/G)[n]\longrightarrow 0,\]
where $\Pic(H/G)\simeq \Hom_k(G(\overline{k}),\Q/\Z(1))$. Thanks to \cite[Theorem 6.8.4, Corollary 6.8.5]{profinitegroups}, we may view a class $\alpha\in\rmH^2(\pi_1^\et(H/G,\overline{1}),\mmu_n)$ as represented by an extension of profinite groups: 
    \[(E)~~:~~ 1\longrightarrow \mmu_n\longrightarrow E\longrightarrow \pi_1^\et(H/G,\overline{1})\longrightarrow 1,\]
     where $\mmu_n$ is viewed as a discrete group, $E$ is profinite, and we require the morphisms to be continuous. Pulling $(E)$ back along the inclusion $G(\overline{k})\hookrightarrow \pi_1^\et(H/G,\overline{1})$ then yields a commutative diagram with exact rows:
\[
    \begin{tikzcd} 
        1 \ar[r] & \mmu_n \ar[r] \ar[d, equal] & \widetilde{E} \ar[d, hookrightarrow] \ar[r] & G(\overline{k}) \ar[d, hookrightarrow] \ar[r]& 1 \\
        1 \ar[r] & \mmu_n \ar[r]& E \ar[r]  & \pi_1^\et(H/G,\overline{1}) \ar[r] & 1,
    \end{tikzcd}
\]
where $\widetilde{E}:=\ker[E\to\Gamma_k]$ comes naturally equipped with an action of $\Gamma_k$ such that the top row is a central extension of $\Gamma_k$-groups (indeed $G(\overline{k})\hookrightarrow \pi_1^\et(H/G,\overline{1})$ acts trivially on $\Q/\Z(1)$ by assumption). In particular, it gives rise to a well-defined morphism $\phi:\rmH^2(\pi_1^\et(H/G,\overline{1}),\mmu_n)\to \mathrm{Ext}_k^c(G(\overline{k}),\mmu_n)$. Consider the section $\delta_1:\Gamma_k\to \pi_1^\et(H/G,\overline{1})$ associated to $1\in (H/G)(k)$, which provides an isomorphism $\pi_1^\et(H/G,\overline{1})\simeq G(\overline{k})\rtimes\Gamma_k$ that we will fix once and for all. We claim that $\phi$ admits a natural section $\mathrm{Ext}_k^c(G(\overline{k}),\mmu_n)\to N\rmH^2(\pi_1^\et(H/G,\overline{1}),\mmu_n)$, where \[N\rmH^2(\pi_1^\et(H/G,\overline{1}),\mmu_n):=\ker\bigg[\delta_1^*:\rmH^2(\pi_1^\et(H/G,\overline{1}),\mmu_n)\to \Br(k)[n]\bigg].\] 
Starting from (the class of) a central extension $(\widetilde{E})$ in $\mathrm{Ext}_k^c(G(\overline{k}),\mmu_n)$,
we may construct an element $\alpha\in N\rmH^2(\pi_1^\et(H/G,\overline{1}),\mmu_n)$ by considering the extension of profinite groups:
\[1\longrightarrow \mmu_n\longrightarrow \widetilde{E}\rtimes \Gamma_k\longrightarrow G(\overline{k})\rtimes \Gamma_k\longrightarrow 1,\]
where the semidirect product structures on both $\widetilde{E}\rtimes \Gamma_k$ and $G(\overline{k})\rtimes \Gamma_k$ are given by the $\Gamma_k$-actions prescribed on $\widetilde{E}$ and $G(\overline{k})$ respectively. Note that the latter action coincides with the conjugation action induced by $\delta_1$. Moreover the map $\widetilde{E}\rtimes \Gamma_k\to G(\overline{k})\rtimes \Gamma_k$ is defined as a morphism of semidirect products in the obvious way, since $\widetilde{E}\to G(\overline{k})$ is supposed to be $\Gamma_k$-equivariant. The isomorphism $\pi_1^\et(H/G,\overline{1})\simeq G(\overline{k})\rtimes\Gamma_k$ allows us to view the class of the above extension as an element of $N\rmH^2(\pi_1^\et(H/G,\overline{1}),\mmu_n)$ (it is clear that this extension pulls back trivially along $\delta_1:\Gamma_k\to \pi_1^\et(H/G,\overline{1})$), so that we obtain a morphism $\mathrm{Ext}_k^c(G(\overline{k}),\mmu_n)\to N\rmH^2(\pi_1^\et(H/G,\overline{1}),\mmu_n)$. Actually, any class $\alpha\in N\rmH^2(\pi_1^\et(H/G,\overline{1}),\mmu_n)$ arises from this construction. Indeed, if we represent $\alpha$ by a central extension $(E)$, then pulling back along $\delta_1$ yields a commutative diagram with exact rows:
    \[
    \begin{tikzcd} 
        1 \ar[r] & \mmu_n \ar[r] \ar[d, equal] & \mmu_n\rtimes \Gamma_k \ar[d] \ar[dr, phantom] \ar[r] & \Gamma_k \ar[d, "\delta_1"] \ar[r]& 1 \\
        1 \ar[r] & \mmu_n \ar[r]& E \ar[r] & \pi_1^\et(H/G,\overline{1}) \ar[r]& 1.
    \end{tikzcd}
\]
In particular we obtain an injection $\Gamma_k\hookrightarrow E$ which, by functoriality of the pullback, is a section of the projection $E\to \pi_1^\et(H/G,\overline{1})\to\Gamma_k$, so that we may write $E=\widetilde{E}\rtimes \Gamma_k$ where $\widetilde{E}:=\ker[E\to \Gamma_k]$. The class of $(\widetilde{E})$ seen as an element of $\mathrm{Ext}_k^c(G(\overline{k}),\mmu_n)$ maps to $\alpha$ along the morphism $\mathrm{Ext}_k^c(G(\overline{k}),\mmu_n)\to N\rmH^2(\pi_1^\et(H/G,\overline{1}),\mmu_n)$ that we constructed. Note that we have a natural commutative diagram with exact rows:
   \[
    \begin{tikzcd} 
        0\ar[r] & \Hom_k(G(\overline{k}),\Q/\Z(1))/n \ar[r] \ar[d, "\sim"] & N\rmH^2(\pi_1^\et(H/G,\overline{1}),\mmu_n) \ar[d, hook] \ar[dr, phantom] \ar[r] & \Br^0(H/G)[n] \ar[d, hook] \ar[r]& 0 \\
        0 \ar[r] & \Pic(H/G)/n \ar[r]& \rmH^2(\pi_1^\et(H/G,\overline{1}),\mmu_n) \ar[r] & \Br(H/G)[n] \ar[r]& 0,
    \end{tikzcd}
\]
where we can now naturally identify $N\rmH^2(\pi_1^\et(H/G,\overline{1}),\mmu_n)$ with $\mathrm{Ext}_k^c(G(\overline{k}),\mmu_n)$. Taking the colimit, over $n\ge 1$, of the collection of diagrams as above, and using the fact that $\Pic(H/G)$ is torsion, we therefore obtain the desired isomorphism $\mathrm{Ext}_k^c(G(\overline{k}),\Q/\Z(1))\simeq \Br^0(H/G)$, which is by construction compatible with the inclusions $\Br^0(H/G)\hookrightarrow \Br(H/G)$ and $\mathrm{Ext}_k^c(G(\overline{k}),\Q/\Z(1))\hookrightarrow\rmH^2(\pi_1^\et(H/G,\overline{1}),\Q/\Z(1))$. The verification of the functoriality in the $k$-group $G$ is straightforward.
    \end{proof}
\subsection{The criterion} We now want to apply the general criterion provided by Theorem~\ref{criterenonramifiehomogene} to this concrete description of the normalised Brauer group. Indeed, if $G$ is a finite $k$-group, let us fix an embedding $G\hookrightarrow \SL_{n,k}$ for some $n\ge 1$. One can then consider the homogeneous space $\SL_{n,k}/G$. In this case we may naturally identify $\Br(\SL_{n,k}/G)=\Br^0(\SL_{n,k}/G)\oplus \Br(k)$ and $\Br^0(\SL_{n,k}/G)\simeq\mathrm{Ext}_k^c(G(\overline{k}),\Q/\Z(1))$, and translate Theorem \ref{criterenonramifiehomogene} in terms of $\Gamma_k$-equivariant central extensions of the $\Gamma_k$-group $G(\overline{k})$ by the Galois module $\Q/\Z(1)$. 

From now on we adopt the following notation: 
\[ \Sha^2_{\mathrm{bic},0}\mathrm{Ext}^c_k(G(\overline{k}),\Q/\Z(1)):= \Sha^2_{\mathrm{bic},0}(G(\overline{k})\rtimes \Gamma_k,\Q/\Z(1))\cap \Br^0(\SL_{n,k}/G)\]
and 
\[\Sha^2_{\mathrm{cyc},\mathrm{pcyc}}\mathrm{Ext}^c_k (G(\overline{k}),\Q/\Z(1)):=\Sha^2_{\mathrm{cyc},\mathrm{pcyc}}(G(\overline{k})\rtimes \Gamma_k,\Q/\Z(1))\cap \Br^0(\SL_{n,k}/G).\]
\begin{thm}\label{nonramificationextension}
Let $G$ be a finite $k$-group together with an embedding of $k$-groups $G\hookrightarrow \SL_{n,k}$ for some $n\ge1$. Let $\alpha\in\Br^0(\SL_{n,k}/G)\simeq \mathrm{Ext}^c_k(G(\overline{k}),\Q/\Z(1))$ be a class represented by a central extension of $\Gamma_k$-groups
\[
  (E)~~: ~~ 1\longrightarrow \Q/\Z(1)\longrightarrow E(\overline{k})\longrightarrow G(\overline{k})\longrightarrow 1.
\]
Then $\alpha$ is unramified if and only if it satisfies the following two conditions:
\begin{enumerate}
    \item[(i)] \emph{(Bogomolov condition).} The extension $(E)$ splits (as an extension of abstract groups) when pulled back to any bicyclic subgroup of $G(\overline{k})$;
    \item[(ii)] \emph{(Galois condition).} For any $\sigma\in\Gamma_k\backslash\{1\}$ and for any two elements $\tau,\gamma\in G(\overline{k})$ such that $\gamma\mkern2.5mu {}^\sigma \mkern-2.5mu\tau\gamma^{-1}=\tau^{\chi(\sigma)}$, there exists a lift of abstract groups $\psi:\langle \tau\rangle\to E(\overline{k})$ and an element $e_\gamma\in E(\overline{k})$ which maps to $\gamma$ and verifies $e_\gamma\mkern1mu {}^\sigma \mkern-1mu\psi(\tau)e_\gamma^{-1}=\psi(\tau)^{\chi(\sigma)}$.
\end{enumerate}
\end{thm}

\begin{rmk}\label{formulealgebrique}
Consider the (normalised) algebraic subgroup \[\Br_1^0(\SL_{n,k}/G)\simeq \Br_1(\SL_{n,k}/G)/\Br(k)\simeq \rmH^1(k,\widehat{G}^{\mathrm{ab}}),\]
which we naturally identify thanks to Proposition \ref{normaliseextension} with the kernel of the base change morphism $\mathrm{Ext}_k^c(G(\overline{k}),\Q/\Z(1))\to \rmH^2(G(\overline{k}),\Q/\Z)$ induced by $\SL_{n,\overline{k}}/G\to \SL_{n,k}/G$. Let us briefly explain how to recover the formula for the algebraic unramified Brauer group discussed in \cite[§3]{harpazwittenbergmassey} and based on the work of Demarche, Harari and Lucchini Arteche (note that the homogeneous spaces considered in \emph{loc.\ cit.}\ do not necessarily admit a $k$-point, however the authors explain why this consideration is not restrictive). Let $\sigma \in\Gamma_k$ and $\tau\in G(\overline{k})$ satisfy the equality ${^\sigma[\tau]} = [\tau^{\chi(\sigma)}]\in G(\overline{k})/(\text{conjugacy})$. Then the
image $\overline{\tau}$ of $\tau$ in $G(\overline{k})^\mathrm{ab}$ satisfies ${^\sigma
\overline{\tau}}=\chi(\sigma)\overline{\tau}$ and is therefore $\sigma$-invariant when viewed as an element of the twisted Cartier dual $\Hom(\widehat{G}^{\mathrm{ab}},\Q/\Z)$ of $\widehat{G}^{\mathrm{ab}}=\Hom(G(\overline{k}),\Q/\Z(1))$ (here $\Q/\Z$ denotes $\Q/\Z(1)$ equipped with the trivial $\Gamma_k$-action). Note that this bidual is canonically isomorphic to $G^{\mathrm{ab}}$ as an abstract group, but the Galois action is twisted by the cyclotomic character $\chi$. If we thus let $\widehat{\Z}$ act on $\Hom(\widehat{G}^{\mathrm{ab}},\Q/\Z)$ via $\sigma$, then we may view $\overline{\tau}$ as an element of $\rmH^0(\widehat{\Z},\Hom(\widehat{G}^{\mathrm{ab}},\Q/\Z))$, and we obtain a composite map
\[\begin{tikzcd}
\rmH^1(k,\widehat{G}^{\mathrm{ab}}) \ar[r, "\sigma^*"] & \rmH^1(\widehat{\Z},\widehat{G}^{\mathrm{ab}})  \ar[r, "\overline{\tau}\smile(-)"] & \rmH^1(\widehat{\Z},\Q/\Z) = \Q/\Z, \end{tikzcd}\]
where $\sigma$ is viewed as map $\widehat{\Z}\to \Gamma_k$ sending $1\mapsto \sigma$.
Taking these identifications into consideration, we may represent a class $\alpha\in\Br_1^0(\SL_{n,k}/G)$ by a central extension: 
\[(E)~~:~~1\longrightarrow\Q/\Z(1)\longrightarrow E(\overline{k})\longrightarrow G(\overline{k})\longrightarrow1\]
where $E(\overline{k})$ is isomorphic to $\Q/\Z(1)\times G(\overline{k})$ as an abstract group, and it is endowed with a $\Gamma_k$-action given by some $1$-cocycle $f_{(-)}\in Z^1(k,\widehat{G}^\mathrm{ab})$ whose class in $\rmH^1(k,\widehat{G}^\mathrm{ab})$ is $\alpha$. The condition ${^\sigma[\tau]}=[\tau^{\chi(\sigma)}]$ in $G(\overline{k})/(\text{conjugacy})$ is equivalent to the existence of some $\gamma\in G(\overline{k})$ such that $\gamma\mkern2.5mu {}^\sigma \mkern-2.5mu\tau\gamma^{-1}=\tau^{\chi(\sigma)}$. By Theorem~\ref{nonramificationextension}, we see that $\alpha$ is unramified if and only if, for each such $\tau,\gamma\in G(\overline{k})$ and $\sigma\in \Gamma_k$, there exists a lift $\psi:\langle\tau\rangle\to E(\overline{k})$ and some $e_\gamma\in E(\overline{k})$ mapping to $\gamma$ and such that $e_\gamma{^\sigma\psi(\tau)}e_\gamma^{-1}=\psi(\tau)^{\chi(\sigma)}$ in $E(\overline{k})$. If we write $\psi(\tau)=(\zeta,\tau)$ and $e_\gamma=(\lambda,\gamma)$ for some $\zeta,\lambda\in \Q/\Z(1)$, then the previous equality amounts to $({^\sigma \zeta}+f_\sigma(\tau),\gamma\mkern2.5mu {}^\sigma \mkern-2.5mu\tau\gamma^{-1})=(\zeta^{\chi(\sigma)},\tau^{\chi(\sigma)})$, which in turn is equivalent to ${^\sigma\zeta}+f_\sigma(\tau)=\zeta^{\chi(\sigma)}$. But ${^\sigma\zeta}=\zeta^{\chi(\sigma)}$ since the $\Gamma_k$-action on $E(\overline{k})$ restricts to the usual action on $\Q/\Z(1)$ via the cyclotomic character. This forces $f_\sigma(\tau)=0$, and \emph{a fortiori} one has 
\[\overline{\tau}\smile \sigma^*\alpha=\overline{\tau}\smile \sigma^*[f_{(-)}]=f_\sigma(\overline{\tau})=0\in \Q/\Z.\] 
This readily yields the formula provided in \cite[Proposition 3.3]{harpazwittenbergmassey}:
\[\begin{aligned} \Br_{\mathrm{nr},1}^0(\SL_{n,k}/G) =  \bigg\{ \alpha\in \rmH^1(k,\widehat{G}^{\mathrm{ab}})~\bigg|  ~\forall \sigma&\in\Gamma_k,~\forall \tau\in G(\overline{k})~\text{such that } {^\sigma[\tau]}=[\tau^{\chi(\sigma)}]\\ & \text{in }G(\overline{k})/(\text{conjugacy}),\text{ then }\overline{\tau}\smile \sigma^*\alpha=0\in\Q/\Z\bigg\}.
\end{aligned}\]
\end{rmk}
\begin{rmk} We would like to highlight the fact that the cohomological formula from Theorem~\ref{criterenonramifiehomogene} for the unramified Brauer group still holds in characteristic $p>0$, as soon as we avoid $p$-torsion. Indeed, the criterion from Proposition \ref{WittenbergBrnr} is valid over any field but the identification, for a field $L$ of dimension $\le 1$ and characteristic $p>0$, between $\Gamma_{L((t))}$ and $ \widehat{\Z}(1)\rtimes\Gamma_L$, does not hold because of Artin-Schreier theory. It does, however, as soon as we pass to a pro-$\ell$-Sylow subgroup, for any prime number $\ell\neq p$. Up to restricting to $\ell$-primary torsion classes in the Brauer group for such $\ell$'s and using a restriction-corestriction argument, the proof of Theorem \ref{criterenonramifiehomogene} works \emph{mutatis mutandis}. Consequently, the proof of Theorem \ref{nonramificationextension} (and hence its conclusion) remains unchanged as soon as we avoid $p$-torsion as well.
    
\end{rmk}
\begin{proof}[Proof of Theorem \ref{nonramificationextension}] Let us consider $\alpha\in \Br^0(\SL_{n,k}/G)$ and write the corresponding central extension:
\begin{equation}\label{extensioncentrale}1\longrightarrow \Q/\Z(1)\longrightarrow E(\overline{k})\overset{\pi}{\longrightarrow} G(\overline{k})\longrightarrow 1.\end{equation}
As in the proof of Proposition \ref{normaliseextension}, we fix the section $\Gamma_k\to \pi_1^\et(\SL_{n,k}/G,\overline{1})$ induced by the neutral element $1\in \SL_{n,k}(k)$ and we view $\alpha$ as an element of $\Br(\SL_{n,k}/G)\simeq \rmH^2(G(\overline{k})\rtimes \Gamma_k,\Q/\Z(1))$. By Corollary \ref{unramifiednormalisehomogene}, we have the equality:
\[\Brnr^0(\SL_{n,k}/G)=\Sha^2_{\mathrm{bic},0}\mathrm{Ext}^c_k(G(\overline{k}),\Q/\Z(1))\cap \Sha^2_{\mathrm{cyc},\mathrm{pcyc}}\mathrm{Ext}^c_k (G(\overline{k}),\Q/\Z(1)).\]
We thus need to show that the following two statements are equivalent:
\begin{enumerate}
\item[(a)] conditions (i) and (ii) from Theorem \ref{nonramificationextension} both hold for the extension \eqref{extensioncentrale};
\item[(b)] the following two conditions are satisfied: 
\begin{enumerate}
    \item[(i')] condition (a), (i) holds for the extension \eqref{extensioncentrale};
    \item[(ii')] for any group $D$ which is an extension
\[1\longrightarrow C\longrightarrow D\longrightarrow\widehat{\Z}\longrightarrow 1\] 
with $C$ a cyclic group, together with a morphism $D\to G(\overline{k})\rtimes\Gamma_k$ which restricts to an injection $C\hookrightarrow G(\overline{k})$, the class of the extension \eqref{extensioncentrale} in $\mathrm{Ext}^c_k(G(\overline{k}),\Q/\Z(1))$ vanishes along the pullback $\rmH^2(G(\overline{k})\rtimes\Gamma_k,\Q/\Z(1))\to \rmH^2(D,\Q/\Z(1))$. 
\end{enumerate}
\end{enumerate}

Actually, we may further assume in (b), (ii') that $\widehat{\Z}$ acts on $C=\langle \tau\rangle$ via the cyclotomic character evaluated on the image $\langle\sigma\rangle$ of $\widehat{\Z}\to\Gamma_k$, i.e.\ ${^1 \tau}=\tau^{\chi(\sigma)}\in C$. The necessity holds trivially thanks to the criterion of Theorem \ref{criterenonramifiehomogene}; for the sufficiency, it follows from the first part of the proof of the latter. Indeed, if $\alpha\in\Br(\SL_{n,k}/G)\simeq \rmH^2(G(\overline{k})\rtimes\Gamma_k,\Q/\Z(1))$ is ramified, then there exists a morphism of semidirect products $\varphi:\widehat{\Z}(1)\rtimes\widehat{\Z}\to G(\overline{k})\rtimes\Gamma_k$ where $\widehat{\Z}$ acts on $\widehat{\Z}(1)$ through the cyclotomic character, and such that $\varphi^*(\alpha)\neq0\in\rmH^2(\widehat{\Z}(1)\rtimes\widehat{\Z},\Q/\Z(1))$. If we let $\zeta\in \widehat{\Z}(1)$ be a topological generator, $(\tau,1):=\varphi(\zeta,0)$ and $(\gamma,\sigma):=\varphi(0,1)$, we can therefore factor $\varphi$ as:
\[\begin{tikzcd}
	{\widehat{\Z}(1)\rtimes\widehat{\Z}} & {\langle\tau\rangle\rtimes \widehat{\Z}} & { G(\overline{k})\rtimes\Gamma_k},
	\arrow["{\varphi_1}", from=1-1, to=1-2]
	\arrow["{\varphi_2}", from=1-2, to=1-3]
\end{tikzcd}\]
by putting $\varphi_1(\zeta,0)=(\tau,0)$, $\varphi_1(0,1)=(1,1)$, $\varphi_2(\tau,0)=(\tau,1)$ and $\varphi_2(1,1)=(\gamma,\sigma)$, and by requesting that $\widehat{\Z}$ acts on $\tau$ via the cyclotomic character evaluated on its image $\langle\sigma\rangle\subset \Gamma_k$. In particular if $\varphi_2^*(\alpha)=0$ for every choice of $\varphi_2$ as above, then $\alpha$ is unramified.

Let us now suppose that \eqref{extensioncentrale} satisfies (b). Then (a) is easily seen to hold. Indeed (i) is satisfied, and if $\sigma\in\Gamma_k\setminus\{1\}$ and  $\tau,\gamma\in G(\overline{k})$ are such that $\gamma\mkern2.5mu {}^\sigma \mkern-2.5mu\tau\gamma^{-1}=\tau^{\chi(\sigma)}$, then we let $D$ be the semidirect product $\langle \tau\rangle\rtimes \widehat{\Z}$ where $\widehat{\Z}$ acts on $\langle\tau\rangle$ via ${^1\tau}=\tau^{\chi(\sigma)}$. By mimicking the proof of Proposition \ref{normaliseextension} and using the fact that $\rmH^2(\widehat{\Z},\Q/\Z(1))=0$, we see that there are natural isomorphisms \[\rmH^2(D,\Q/\Z(1))\simeq N\rmH^2(D,\Q/\Z(1))\simeq \mathrm{Ext}_{\widehat{\Z}}^c(C,\Q/\Z(1)),\]
where we recall that the latter group denotes classes of central extensions of $\widehat{\Z}$-groups, where $\widehat{\Z}$ acts on both $C$ and $\Q/\Z(1)$ via the cyclotomic character. We define a morphism $\varphi:D\to G(\overline{k})\rtimes\Gamma_k$ by letting $\varphi(\tau,0)=(\tau,1)$ and $\varphi(0,1)=(\gamma,\sigma)$. It is easy to see that $\varphi$ is a morphism of semidirect products, and hence we may apply condition (ii'). The pullback map $\phi^*:\rmH^2(G(\overline{k})\rtimes\Gamma_k,\Q/\Z(1))\to \rmH^2(D,\Q/\Z(1))$ restricts to a map \[\phi^*:\mathrm{Ext}_k^c(G(\overline{k}),\Q/\Z(1))\to \mathrm{Ext}^c_{\widehat{\Z}}(C,\Q/\Z(1))\] 
as follows: the morphism $D\to G(\overline{k})\rtimes \Gamma_k$ induces a group morphism $\phi:C\to G(\overline{k})$ which is $\gamma$-twisted-equivariant, in the sense that it satisfies $\phi({^1\tau})=\gamma\mkern2.5mu {}^\sigma \mkern-2.5mu\phi(\tau)\gamma^{-1}$. The map $\varphi^*$ then consists of pulling back the central extension \eqref{extensioncentrale} along this morphism; in other words, we have a commutative diagram of extensions
    \[
    \begin{tikzcd} 
        1 \ar[r] & \Q/\Z(1) \ar[r] \ar[d, equal] & E_C \ar[d, "\phi_E"] \ar[dr, phantom] \ar[r] & C \ar[d, "\phi"] \ar[r]& 1 \\
        1\ar[r] & \Q/\Z(1) \ar[r]& E(\overline{k}) \ar[r] & G(\overline{k}) \ar[r]& 1,
    \end{tikzcd}
\]
where the groups involved in the top row are equipped with their $\widehat{\Z}$-action, those of the bottom one are equipped with their $\Gamma_k$-action, such that there exists some $e_\gamma\in E(\overline{k})$ for which $\pi(e_\gamma)=\gamma$ and such that for any $e\in E_C$, we have $\phi_E({^1 e})=e_\gamma\mkern2.5mu{}^\sigma \mkern-2.5mu\phi_E(e)e_\gamma^{-1}$. Since the pullback extension in the above diagram splits $\widehat{\Z}$-equivariantly by assumption, there exist an element $e_\tau$ of $E(\overline{k})$ that maps to $\tau$ in $G(\overline{k})$ with the same order, and subject to the semi-commutation condition: 
\[e_\gamma\mkern2.5mu {}^\sigma \mkern-2.5mu e_\tau e_\gamma^{-1}=e_\tau^{\chi(\sigma)}\in E(\overline{k}).\] The morphism $\psi:\langle\tau\rangle\to E(\overline{k})$ given by sending $\tau\mapsto e_\tau$ then yields the desired lift in condition (a), (ii).

Conversely, assume that \eqref{extensioncentrale} satisfies (a) and let $\phi:D\to G(\overline{k})\rtimes\Gamma_k$ be the prescribed morphism. Of course, (i') is tautologically satisfied. On the other hand, let $D$ be as described in condition (ii'), let $\tau\in C$ be a generator and $\sigma\in \Gamma_k$ be the image of $1\in\widehat{\Z}$ via $D\to G(\overline{k})\rtimes \Gamma_k\to \Gamma_k$. Fix any section $\widehat{\Z}\to D$ (which exists as $\widehat{\Z}$ is a free profinite group), write $D=C\rtimes\widehat{\Z}$ and denote by $(\gamma,\sigma)$ the image of $(0,1)\in D$ in $G(\overline{k})\rtimes\Gamma_k$. We obtain that $\mathrm{Im}(\phi)\subset G(\overline{k})\rtimes\Gamma_k$ is generated by $(\tau,1)$ and $(\gamma,\sigma)$, subject to the relation
\[(\gamma,\sigma)(\tau,1)(\gamma,\sigma)^{-1}=(\tau^{\chi(\sigma)},1).\]
By condition (a), (ii), we know that there exists an abstract lift $\psi:\langle\tau\rangle\to E(\overline{k})$ and an element $e_\gamma\in E(\overline{k})$ which maps to $\gamma\in G(\overline{k})$ and satisfies $e_\gamma\mkern1mu {}^\sigma \mkern-1mu\psi(\tau)e_\gamma^{-1}=\psi(\tau)^{\chi(\sigma)}$. In particular, we may construct a $\gamma$-twisted-equivariant lift $C\to E(\overline{k})$ of the map $\phi:C\to G(\overline{k})$ in the natural way, by sending $\tau\mapsto\psi(\tau)$. This implies that the class of \eqref{extensioncentrale} is zero in $\mathrm{Ext}_{\widehat{\Z}}^c(C,\Q/\Z(1))\simeq\rmH^2(D,\Q/\Z(1))$, so that (ii') is satisfied and (b) holds.
\end{proof}

\begin{cor}\label{nonramificationextensionconstant}
Suppose that $G$ is a constant $k$-group, let $N$ be its exponent, and assume that $\mmu_N\subset k$. Let $\alpha\in\Br^0(\SL_{n,k}/G)\simeq \mathrm{Ext}^c_k(G,\Q/\Z(1))$ be a class represented by a central extension of $\Gamma_k$-groups
\[
  (E)~~:~~ 1\longrightarrow \Q/\Z(1)\longrightarrow E(\overline{k})\longrightarrow G\longrightarrow 1.
\]
Then $\alpha$ is unramified if and only if it satisfies the following two conditions:
\begin{enumerate}
    \item[(i)] \emph{(Bogomolov condition).} The extension $(E)$ splits (as an extension of abstract groups) when pulled back to any bicyclic subgroup of $G$ ;
    \item[(ii)] \emph{(Galois condition).} For any $\sigma\in\Gamma_k\backslash\{1\}$ and for any cyclic subgroup $C\subset G$, the extension $(E)$ splits $\sigma$-equivariantly when pulled back to $C$.
\end{enumerate}
\end{cor}
\begin{proof}
We need to show the equivalence between the following two statements:
\begin{enumerate}
\item[(a)] conditions (i) and (ii) from Theorem \ref{nonramificationextension} both hold for the extension (E);
\item[(b)] the following two conditions are satisfied: 
\begin{enumerate}
    \item[(i')] condition (a), (i) holds for the extension (E);
    \item[(ii')] For any $\sigma\in\Gamma_k\backslash\{1\}$ and for any cyclic subgroup $C\subset G$, the extension $(E)$ splits $\sigma$-equivariantly when pulled back to $C$.
\end{enumerate}
\end{enumerate}

Note that in the relations of the form $\gamma{^\sigma \tau}\gamma^{-1}=\tau^{\chi(\sigma)}$ from condition (a), (ii), the exponent on the right hand side only depends on the modulo $N$ cyclotomic character $\Gamma_k\to\mathrm{Aut}(\mmu_N)$. If we thus assume that $\mmu_N\subset k$, then this character is trivial and we have the equality $\gamma\tau\gamma^{-1}=\gamma\mkern2.5mu {}^\sigma \mkern-2.5mu\tau\gamma^{-1}=\tau$, i.e.\ the subgroup $\langle \tau,\gamma\rangle\subset G$ is bicyclic. 

Assume (b). If $\tau,\gamma\in G$ are such that $\gamma\tau\gamma^{-1}=\tau$, then condition (b), (i') gives an abstract lift $\phi:\langle\tau,\gamma\rangle\to E(\overline{k})$. Let $e_\tau,e_\gamma\in E(\overline{k})$ be two (\emph{a priori} different) elements mapping to $\tau$ and $\gamma$ respectively. Since $e_\tau$ and $\phi(\tau)$, resp.\ $e_\gamma$ and $\phi(\gamma)$ differ by a central element, we obtain that $e_\tau$ and $e_\gamma$ commute. On the other hand, (ii') provides a $\sigma$-equivariant lift $\psi:\langle\tau\rangle\to E(\overline{k})$. We just need to check that we may choose some $e_\gamma\in E(\overline{k})$ satisfying $e_\gamma\mkern1mu {}^\sigma \mkern-1mu\psi(\tau)e_\gamma^{-1}=\psi(\tau)$. Since $e_\gamma$ and $\psi(\tau)$ commute, this amounts to asking whether ${^\sigma \psi(\tau)}=\psi(\tau)$, which is true since $\psi$ is $\sigma$-equivariant. This implies (a), (ii).  

Conversely, since $\alpha$ is unramified, then it satisfies (a) thanks to Theorem \ref{nonramificationextension}. Let $\sigma\in\Gamma_k\setminus\{1\}$ and $\tau\in G$. We may let $\gamma=1$ in condition (a), (ii) and we thus obtain a lift $\psi:\langle\tau\rangle\to E(\overline{k})$ and some element $e\in \Q/\Z(1)\subset E(\overline{k})$ which satisfies $e\mkern1mu {}^\sigma \mkern-1mu\psi(\tau)e^{-1}=\psi(\tau)$. Since $e$ is central, we thus have ${^\sigma\psi(\tau)}=\psi(\tau)$. In other words $\psi$ is $\sigma$-equivariant, and we obtain condition (b), (ii'). Since (i') is tautological, we obtain (b). This concludes the proof.
\end{proof}

Let us briefly take a look at the algebraic subgroup of the unramified Brauer group under the considerations of Theorem \ref{nonramificationextensionconstant}. We can recover the following vanishing result due to Lucchini Arteche \cite[Proposition 5.9]{lucchiniartechenonramifie}:
\begin{cor}\label{nonramificationextensionconstantalgebrique}
    Under the hypotheses of Corollary \ref{nonramificationextensionconstant}, we have $\Br_{\mathrm{nr},1}^0(\SL_{n,k}/G)=0$.
\end{cor}
\begin{proof} An element of $\Br_1^0(\SL_{n,k}/G)$ is represented by the class of a central extension 
\[1\longrightarrow \Q/\Z(1)\longrightarrow E(\overline{k})\longrightarrow G\longrightarrow 1\]
where $E(\overline{k})$ is isomorphic to $\Q/\Z(1)\times G$ as an abstract group. A $\Gamma_k$-action on $E(\overline{k})$ which makes the arrows in the above extension $\Gamma_k$-equivariant therefore corresponds to a function $f_{(-)}:\Gamma_k\to \Hom(G,\Q/\Z(1))\simeq \Hom(G,\mmu_N)$ which is a $1$-cocycle; since $\Gamma_k$ acts trivially on $\Hom(G,\mmu_N)$, this cocycle is just a morphism, hence defines an element of $\Hom(\Gamma_k,\Hom(G,\mmu_N))$. Fix $\tau \in G$ and $\sigma\in \Gamma_k\setminus\{1\}$; if there exists a $\sigma$-equivariant section $\psi:\langle\tau\rangle\to E(\overline{k})$, we may write $\psi(\tau)=(\lambda,\tau)$ for some $\lambda\in \mmu_N$. We then have ${^\sigma(\lambda,\tau)}=(\lambda,\tau)$; but ${^\sigma(\lambda,\tau)}=(\lambda+f_\sigma(\tau),\tau)$ since $\Gamma_k$ acts trivially on $\mmu_N$ and on $G$, so this forces $f_\sigma(\tau)=0$. Applying the same argument for varying choices of $\tau\in G$ and $\sigma\in \Gamma_k$, we obtain that $f_{(-)}$ is identically zero, i.e.\ $E(\overline{k})\simeq \Q/\Z(1)\times G$ as a $\Gamma_k$-group. Its class in $\Br_{\mathrm{nr},1}^0(\SL_{n,k}/G)$ is therefore trivial. 
\end{proof}
\section{Algorithms for computing the unramified Brauer group} 

Let now $G$ be a finite $k$-group and let $G\hookrightarrow \SL_{n,k}$ be an embedding of $k$-groups for some $n\ge 1$. A wide variety of sources in the literature describe an algorithm, or at least an effective formula in order to compute the algebraic unramified Brauer group of $\SL_{n,k}/G$, see e.g.\ \cite{lucchiniartechenonramifie} or \cite{borovoidemarcheharari}. On the other hand, the computation of the Bogomolov multiplier of $G(\overline{k})$ has been investigated by Bogomolov in his original paper \cite{bogomolov}, and more recently in e.g.\ \cite{KreschTschinkel24} and \cite{Moravec12}. However these papers do not discuss the existence of an algorithm for the computation of the full unramified Brauer group of $\SL_{n,k}/G$, notably its transcendental part when $k$ is not algebraically closed. Here we present a method to compute the unramified Brauer group over an arbitrary field of characteristic zero:
\begin{thm}\label{algorithmenonramifie}
    There exists an algorithm that, given the datum of a field $k$ of characteristic~$0$ and a finite $k$-group $G$, computes the group $\Brnr^0(\SL_{n,k}/G)$ for any embedding $G\hookrightarrow\SL_{n,k}$ for some $n\ge1$.
\end{thm}
\begin{rmk}
    We want to stress the fact that $\Brnr(\SL_{n,k}/G)\simeq \Brnr^0(\SL_{n,k}/G)\oplus \Br(k)$ and $\Brnr^0(\SL_{n,k}/G)$ is a subgroup of $\Brnr(\SL_{n,k}/G)$, and that this algorithm computes the unramified Brauer group in the \emph{strong sense}, that is, it explicitly identifies the elements of $\Brnr(\SL_{n,k}/G)$ as sums of classes of central extensions $\alpha \in\mathrm{Ext}_k^c(G(\overline{k}),\Q/\Z(1))$ and constant classes $\beta\in \Br(k)$. 
\end{rmk}
\subsection{A reduction step} We first show that one can precisely determine all of the unramified normalised Brauer classes of $\SL_{n,k}/G$ as soon as one replaces $k$ by a field with sufficiently many roots of unity over which the $k$-group $G$ splits.
\begin{prop}\label{nonramifieracinesunite}
Let $G$ be a finite constant $k$-group. Fix an embedding $G\hookrightarrow\SL_{n,k}$ for some $n\ge1$ and let $N$ be the order of $G$. Then the natural map:
\[\Brnr^0(\SL_{n,k(\zeta_{N^2})}/G)\longrightarrow B_0(G)\]
is an isomorphism.
\end{prop}
\begin{proof}
Let $K:=k(\zeta_{N^2})$. By Kummer theory \cite[Exposé XVII, Prop. A.2.1, a)]{SGA3} the inclusion $\mmu_N\subset \Q/\Z(1)$ yields a canonical surjection: 
\[\begin{tikzcd}
	{\mathrm{Ext}^c_K(G,\mmu_N)} & {\mathrm{Ext}^c_K(G,\Q/\Z(1))[N]\simeq \Br^0(\SL_{n,K}/G)},
	\arrow[two heads, from=1-1, to=1-2]
\end{tikzcd}\] whose kernel is $Z^1(G,\Q/\Z(1))=\Hom_K(G,\Q/\Z(1))$ since $G$ acts trivially on $\Q/\Z(1)$ (note also that the inclusion $\Hom_K(G,\mmu_N)\hookrightarrow\Hom_K(G,\Q/\Z(1))$ is an isomorphism since $G$ is killed by $N$). The choice of a primitive $(N^2)$\textsuperscript{th}-root of unity $\zeta_{N^2}\in K$ yields isomorphisms of $\Gamma_K$-groups $\mmu_N\simeq \Z/N$ and $\mmu_{N^2}\simeq \Z/N^2$. We then have a commutative diagram:
\[
    \begin{tikzcd} 
          \mathrm{Ext}^c_K(G,\mmu_N) \ar[d, twoheadrightarrow] \ar[r] & \rmH^2(G,\Z/N)\ar[d, twoheadrightarrow]  \\
        \mathrm{Ext}^c_K(G,\Q/\Z(1)) \ar[r] & \rmH^2(G,\Q/\Z)
    \end{tikzcd}
\]
where the horizontal arrows are induced by the base change map $\SL_{n,\overline{K}}/G\to\SL_{n,K}/G$ (after fixing a compatible system of roots of unity in $\overline{k}$), and the surjectivity of the right vertical map comes from the short exact sequence of trivial $G$-modules
\begin{equation}\label{G-modulestriviaux}1\longrightarrow \Z/N\longrightarrow \Q/\Z\overset{\cdot N}{\longrightarrow} \Q/\Z\longrightarrow 1.\end{equation}
The top horizontal map is surjective and admits a section by viewing any extension of abstract groups as an extension of $\Gamma_K$-groups with trivial $\Gamma_K$-action. As for the bottom horizontal arrow, it sends $\Brnr^0(\SL_{n,K}/G)$ to $B_0(G)$ thanks to the functoriality of the unramified Brauer group. On the other hand, Corollary \ref{nonramificationextensionconstantalgebrique} shows that $\Br^0_{\mathrm{nr},1}(\SL_{n,K}/G)$ vanishes, which guarantees the injectivity of $\Brnr^0(\SL_{n,K}/G)\to B_0(G)$. Conversely, let us start with a class $\alpha\in B_0(G)$; we may lift it to some $\widetilde{\alpha}\in \rmH^2(G,\Z/N)$ and, by the above observation, view it as an element of $\mathrm{Ext}^c_K(G,\mmu_N)$. We thus only need to check that the image of $\widetilde{\alpha}$ in $\mathrm{Ext}^c_K(G,\Q/\Z(1))$ is unramified. Note that this class maps to $B_0(G)$ and hence trivially satisfies condition (i) from Corollary \ref{nonramificationextensionconstant}. We claim that it also satisfies condition~(ii). Let $C\subset G$ be a cyclic subgroup and $\sigma\in \Gamma_K\setminus\{1\}$. We may use the functoriality of the Kummer sequence:
\[
    \begin{tikzcd} 
          \mathrm{Ext}^c_K(G,\mmu_N) \ar[d, twoheadrightarrow] \ar[r, "\theta_C^\sigma"] & \mathrm{Ext}^c_{\overline{K}^\sigma}(C,\mmu_N)\ar[d, twoheadrightarrow]  \\
        \mathrm{Ext}^c_K(G,\Q/\Z(1)) \ar[r, "\theta_C^\sigma"] & \mathrm{Ext}^c_{\overline{K}^\sigma}(C,\Q/\Z(1)),
    \end{tikzcd}
\]
where $\theta_C^\sigma$ denotes the natural restriction. Note that the image of $\widetilde{\alpha}$ in $\mathrm{Ext}^c_{\overline{K}^\sigma}(C,\mmu_N)$ factors (by construction) through the inclusion $\rmH^2(C,\Z/N)\hookrightarrow \mathrm{Ext}_K^c(C,\mmu_N)$ that we described earlier. On the other hand, there is a cohomology exact sequence
\[1\longrightarrow \Hom(C,\Q/\Z)/N\longrightarrow \rmH^2(C,\Z/N)\longrightarrow \rmH^2(C,\Q/\Z)\]
associated to the short exact sequence \eqref{G-modulestriviaux}. Since $\Hom(C,\Q/\Z)/N\simeq \Hom(C,\Z/N)$ and $\rmH^2(C,\Q/\Z)$ vanishes, we get that the canonical map $\Hom(C,\Z/N)\to \rmH^2(C,\Z/N)$ is an isomorphism. An element of the latter group is thus obtained by pulling back the extension of abstract groups
\[1\longrightarrow \Z/N\longrightarrow \Z/N^2\longrightarrow \Z/N\longrightarrow 1\]
along a given morphism $C\to \Z/N$. On the other hand, let $\kappa:\Hom_{\overline{K}^\sigma}(C,\mmu_N)\to \mathrm{Ext}^c_{\overline{K}^\sigma}(C,\mmu_N)$ be the map coming from the Kummer sequence. As the image of this map is the kernel of the pushforward $\mathrm{Ext}^c_{\overline{K}^\sigma}(C,\mmu_N)\to \mathrm{Ext}^c_{\overline{K}^\sigma}(C,\Q/\Z(1))$, an extension which lies in this kernel is the pullback, in the category of $\Gamma_K$-groups, of the extension
\[1\longrightarrow \mmu_N\longrightarrow \mmu_{N^2}\longrightarrow \mmu_N\longrightarrow 1\]
along a morphism $C\to \mmu_N$ (where both are viewed as constant $\Gamma_K$-groups), see the proof of \cite[Exposé XVII, Prop. A.2.1, a)]{SGA3}.
Since $\mmu_N$ and $\Z/N$ (resp.\ $\mmu_{N^2}$ and $\Z/N^2$) are isomorphic as $\Gamma_K$-groups, the respective images of $\kappa:\Hom_{\overline{K}^\sigma}(C,\mmu_N)\to \mathrm{Ext}^c_{\overline{K}^\sigma}(C,\mmu_N)$ and the composite map
\[\begin{tikzcd}
	{\rmH^2(G,\Z/N)} & {\mathrm{Ext}^c_K(G,\mmu_N)} & {\mathrm{Ext}^c_{\overline{K}^{\sigma}}(C,\mmu_N)}
	\arrow[hook, from=1-1, to=1-2]
	\arrow[from=1-2, to=1-3, "\theta_C^\sigma"]
\end{tikzcd}\] coincide. Therefore, the image of $\widetilde{\alpha}$ in $\mathrm{Ext}_K^c(G,\Q/\Z(1))$ dies in $\mathrm{Ext}^c_{\overline{K}^\sigma}(C,\Q/\Z(1))$ and hence satisfies condition (ii) from Corollary \ref{nonramificationextensionconstant}, as desired.
\end{proof}
\begin{rmk} This proposition complements the result of Lucchini Arteche \cite[Proposition 5.9]{lucchiniartechenonramifie}. Indeed, if $N=\mathrm{exp}(G)$ and $K=k(\zeta_{N})$ (or any field containing the $N$\textsuperscript{th} roots of unity in $\overline{k}$), then $\Br^0_{\text{nr},1}(\SL_{n,K}/G)$ vanishes. In particular this yields an injection $\Brnr^0(\SL_{n,K}/G)\hookrightarrow \Brnr(\SL_{n,\overline{K}}/G)\simeq B_0(G)$. Proposition~\ref{nonramifieracinesunite} shows that adding slightly more roots of unity to $K$ suffices for this map to be actually onto. One could possibly adapt the proof of Proposition \ref{nonramifieracinesunite} and replace $k(\zeta_{|G|^2})$ by $k(\zeta_{\mathrm{exp}(G)^2})$ in its statement by requiring that the exponent of $B_0(G)$ divide $\mathrm{exp}(G)$; however, this is only conjectured, see \cite{moravecexponent}.
\end{rmk}
\begin{rmk} The assumption that we make on the minimal field $K$ in order to guarantee the validity of the proof of Proposition \ref{nonramifieracinesunite} seems to be optimal. Let for instance $K=\R$ be the field of real numbers, let $\sigma\in\Gamma_\R=\Gal(\C/\R)$ be the complex conjugation, let $G=\Z/2=\{0,1\}$ (viewed as a constant $\Gamma_\R$-group), and consider the extension of constant $\Gamma_\R$-groups
\[1\longrightarrow \mmu_2\longrightarrow \Z/4\longrightarrow \Z/2\longrightarrow 1\]
seen as an element of $\mathrm{Ext}_\R^c(G,\mmu_2)$ (where we write $\mmu_2$ instead of $\Z/2$ to avoid any confusion with the group $G$). This class corresponds to the image of $\mathrm{Id}\in\Hom(G,\Z/2)$ along $\Hom(G,\Z/2)\to \rmH^2(G,\Z/2)\hookrightarrow\mathrm{Ext}_\R^c(G,\Z/2)$, but it is not the image of $\mathrm{Id}$ along the Kummer map $\Hom(G,\Z/2)\to \mathrm{Ext}_\R^c(G,\Z/2)$. Consider the following commutative diagram of extensions of $\Gamma_\R$-groups:
\[\begin{tikzcd}
	1 & {\mmu_2} & {\Z/4} & {\Z/2} & 1 \\
	1 & {\mmu_{\infty}} & {\mmu_\infty\times \Z/2} & {\Z/2} & 1,
	\arrow[from=1-1, to=1-2]
	\arrow[from=1-2, to=1-3]
	\arrow[hook, from=1-2, to=2-2]
	\arrow[from=1-3, to=1-4]
	\arrow[from=1-3, to=2-3]
	\arrow[from=1-4, to=1-5]
	\arrow[Rightarrow, no head, from=1-4, to=2-4]
	\arrow[from=2-1, to=2-2]
	\arrow[from=2-2, to=2-3]
	\arrow[from=2-3, to=2-4]
	\arrow[from=2-4, to=2-5]
\end{tikzcd}\]
where the leftmost vertical map is the natural inclusion and the middle one sends $1\in \Z/4$ to $(i,1)\in \mmu_\infty\times \Z/2$, where $i^2=-1$. The $\Gamma_\R$-action on $\mmu_\infty\times \Z/2$ is defined as ${^{\sigma}(a,b)}:=({(-1)^b}\cdot {^{\sigma}a},b)$ (note that this restricts to complex conjugation on $\mmu_\infty\hookrightarrow \mmu_\infty\times \Z/2$). One can check that all of the arrows in the above diagram are $\Gamma_\R$-equivariant, so that the class of the bottom extension is precisely the image of the class of the top one along the pushforward map $\mathrm{Ext}_\R^c(\Z/2,\mmu_2)\to \mathrm{Ext}_\R^c(\Z/2,\mmu_\infty)$. We claim, however, that the bottom extension does not split $\Gamma_\R$-equivariantly. In fact, an element of the form $(a,1)\in\mmu_\infty\times \Z/2$ of order $2$ satisfies $^\sigma a=a$, so that $^\sigma(a,1)=(-a,1)$, and therefore cannot be $\Gamma_\R$-invariant.
\end{rmk}
\begin{cor}\label{actiontrivialenonramifiée}
Let $N$ be the order of the constant $k$-group $G$ and let $\alpha\in\Br^0(\SL_{n,k}/G)$ be represented by a central extension of $\Gamma_k$-groups
\[(E)~~: ~~1\longrightarrow \Q/\Z(1)\longrightarrow E(\overline{k})\longrightarrow G\longrightarrow 1.\]
If $\alpha$ is unramified, then over $K:=k(\zeta_{N^2})$ the above extension fits into a commutative diagram of $\Gamma_K$-groups with exact rows:
\[ 
    \begin{tikzcd} 
         1 \ar[r] & \mmu_N \ar[d] \ar[r] & \widetilde{E}(\overline{K}) \ar[d]\ar[r] & G \ar[d, equal]\ar[r] & 1 \\
        1 \ar[r] & \Q/\Z(1) \ar[r]& E(\overline{K}) \ar[r] & G \ar[r]& 1.
    \end{tikzcd}
\]
Moreover, for any such diagram, the $\Gamma_K$-group $\widetilde{E}(\overline{K})$ is constant.
\end{cor}
\begin{proof}
The proof of Proposition \ref{nonramifieracinesunite} implies that the groups $\Brnr^0(\SL_{n,K}/G)$ and $B_0(G)$ have isomorphic inverse images in $\mathrm{Ext}_K^c(G,\mmu_N)$ and $\rmH^2(G,\Z/N)$ respectively. In particular any lift $\widetilde{\alpha}\in\mathrm{Ext}_K^c(G,\mmu_N)$ of an unramified class $\alpha\in \Br^0(\SL_{n,K}/G)$ is represented by a constant extension of $\Gamma_K$-groups.
\end{proof}

\subsection{Proof of Theorem \ref{algorithmenonramifie}}

Now, given the datum of a field $k$ and a finite $k$-group $G$ together with an embedding $G\hookrightarrow \SL_{n,k}$, the computation of the unramified normalised Brauer group of $\SL_{n,k}/G$ can be conducted as follows: we will prove the existence of the following sequence of algorithms, which combine to provide the announced algorithm in the statement of the theorem.
\begin{lem}\label{algoA}
    There exists an algorithm that returns, for any finite (abstract) group $G$ of order $N$, a finite list of (isomorphism classes of) central extensions:
\[(\widetilde{E})~~:~~ 1 \longrightarrow \Z/N \longrightarrow \widetilde{E} \longrightarrow G \longrightarrow 1\]
such that this list identifies (as a set) with the inverse image in $\rmH^2(G,\Z/N)$ of $B_0(G)$ along the pushforward map associated to the inclusion $\Z/N\hookrightarrow \Q/\Z$ (in particular the elements of $B_0(G)$ precisely arise as pushforwards of elements of this list).
\end{lem}
\begin{proof} The group $\rmH^2(G,\Z/N)$ is finite and surjects onto $\rmH^2(G,\Q/\Z)$ (since the latter is killed by $N$, this follows from the long exact cohomology sequence associated to the short exact sequence of trivial $G$-modules $1\to \Z/N\to \Q/\Z \to \Q/\Z\to 1$). Describing all the (classes of) extensions of the form
\[(\widetilde{E})~~:~~ 1\longrightarrow \Z/N\longrightarrow \widetilde{E}\longrightarrow G\longrightarrow 1\]
is straightforward: since $\Z/N$ and $G$ are fixed, then the order of $\widetilde{E}$ is fixed. Now we can enumerate all the finite groups of this fixed order, all the injections $\Z/N\to \widetilde{E}$ and surjections $\widetilde{E}\to G$, and check among those which ones fit into an exact sequence. We may then enumerate all of the bicyclic subgroups of the amalgamated sum $\Q/\Z\star_{\Z/N}\widetilde{E}$ whose order is bounded by $|G|$: these groups are finitely many since an element of $\Q/\Z\star_{\Z/N}\widetilde{E}$ whose order is bounded by $|G|$ lifts to the sum of an element $e\in \widetilde{E}$ and an element $\lambda\in \Q/\Z$ which is killed by $|G|\cdot N$). We then check which of these finitely many bicyclic subgroups admit a section. This process eventually stops. (For more effective algorithms, see also \cite[§8]{KreschTschinkel24} (general approach) or \cite[§8]{Moravec12} (in the case of solvable~$G$).)\end{proof}

\begin{lem}\label{algoB}
    There exists an algorithm that, given a field $k$ of characteristic $0$ and a finite $k$-group $G$ of order $N$, returns a finite list of central extensions
\[(\widetilde{E})~~:~~ 1 \longrightarrow \mmu_N \longrightarrow \widetilde{E} \longrightarrow G(\overline{k}) \longrightarrow 1\]
of groups equipped with an action of $\Gal(K/k)$, where $K/k$ is the smallest field extension which splits both $G$ and $\mmu_{N^2}$, such that the arrows are $\Gal(K/k)$-equivariant, and such that:

\begin{enumerate}
\item[(a)] if $\alpha\in\mathrm{Ext}_k^c(G(\overline{k}),\Q/\Z(1))$ denotes the (class of the) pushforward of any such extension along $\mmu_N\hookrightarrow \Q/\Z(1)$, then the image of $\alpha$ along the natural map $\mathrm{Ext}_k^c(G(\overline{k}),\Q/\Z(1))\to \rmH^2(G(\overline{k}),\Q/\Z)$ lies in $B_0(G(\overline{k}))$;
    \item[(b)] the image in $\mathrm{Ext}^c_k(G(\overline{k}),\Q/\Z(1))$ of any element of
$\Br^0_{nr}(\SL_{n,k}/G)$ (via the isomorphism given by Proposition \ref{normaliseextension})
is obtained from one of these extensions by pushing forward along the inclusion $\mmu_N\hookrightarrow\Q/\Z(1)$.
\end{enumerate}
\end{lem}
\begin{proof} Let $E/k$ be the minimal finite Galois extension over which $G$ becomes constant, that is, the kernel of $\Gamma_k\to \mathrm{Aut}(G(\overline{k}))$. We can then take $K$ to be $E(\zeta_{N^2})$. Let $\alpha\in \Brnr^0(\SL_{n,k}/G)\simeq\mathrm{Ext}_k^c(G(\overline{k}),\Q/\Z(1))$. Since $G$ becomes constant over $K$, we know by Corollary \ref{actiontrivialenonramifiée} that the image $\beta$ of $\alpha$ in $\mathrm{Ext}_K^c(G(\overline{K}),\Q/\Z(1))$ must lift to constant central extensions of $\Gamma_K$-groups in $\mathrm{Ext}^c_K(G(\overline{K}),\mmu_N)$. More precisely, if we denote by $(E)_K$ an extension representing $\beta$, then it fits into a commutative diagram of extensions:
\[ 
    \begin{tikzcd} 
    (\widetilde{E})_K~~:~~     1 \ar[r] & \mmu_N \ar[d] \ar[r] & \widetilde{E}(\overline{K}) \ar[d]\ar[r] & G(\overline{K}) \ar[d, equal]\ar[r] & 1 \\
     (E)_K~~:~~   1 \ar[r] & \Q/\Z(1) \ar[r]& E(\overline{K}) \ar[r] & G(\overline{K}) \ar[r]& 1,
    \end{tikzcd}
\]
where $\widetilde{E}(\overline{K})$ is a constant $\Gamma_K$-group. In other words, $(E)_K$ is the pushforward of the central extension of constant $\Gamma_K$-groups $(\widetilde{E})_K$ along the inclusion $\mmu_N\hookrightarrow \Q/\Z(1)$. Note that $(\widetilde{E})_K$ uniquely defines an element of $\rmH^2(G(\overline{K}),\mmu_N)\simeq \rmH^2(G(\overline{K}),\Z/N)$ (since the $\Gamma_K$-action is trivial), and since $(E)_K$ maps to an element of $B_0(G(\overline{K}))$ we may apply Lemma \ref{algoA} to the abstract group $G(\overline{K})$ in order to obtain a list of the finitely many possibilities for $(\widetilde{E})_K$. Let us now move back to the class~$\alpha$. We may also lift it to an element $\widetilde{\alpha}\in\mathrm{Ext}_k^c(G(\overline{k}),\mmu_N)$ which we represent by a central extension~$(\widetilde{E})$. By the functoriality of the Kummer sequence, the image of $\widetilde{\alpha}$ along the base change $\mathrm{Ext}_k^c(G(\overline{k}),\mmu_N)\to \mathrm{Ext}_K^c(G(\overline{K}),\mmu_N)$ is represented by the extension of constant $\Gamma_K$-groups $(\widetilde{E})_K$. This shows that $(\widetilde{E})$ is an extension of $\Gamma_k$-groups for which the action of $\Gamma_k$ is trivial on the subgroup $\Gamma_K\subset \Gamma_k$, and hence factors through $\Gal(K/k)$. We may thus proceed as follows: given a finite list of central extensions of abstract groups: 
\[1\longrightarrow \mmu_N\longrightarrow \widetilde{E}\longrightarrow G(\overline{k})\longrightarrow 1\]
(after amending an isomorphism of groups $\mmu_N\simeq \Z/N$), we may enumerate all the possible actions of $\Gal(K/k)$ on $\widetilde{E}$ such that the maps $\mmu_N\to \widetilde{E}$ and $\widetilde{E}\to G(\overline{k})$ are $\Gal(K/k)$-equivariant (where $\mmu_N$ (resp.\ $G(\overline{k})$) is equipped with the usual action via the cyclotomic character (resp.\ its prescribed $\Gal(K/k)$-action)). This is achieved in a finite number of steps since both $\Gal(K/k)$ and $\widetilde{E}$ are finite, and we obtain the desired list of central extensions that lift $\alpha$ and all of the other unramified Brauer classes.
\end{proof}

\begin{proof}[Proof of Theorem \ref{algorithmenonramifie}] It only remains to show that there exists an algorithm that, given a field $k$ of characteristic $0$ and a finite $k$-group $G$ of order $N$, gives among the finite list of central extensions
\[(E)~~:~~1 \longrightarrow \Q/\Z(1) \longrightarrow E(\overline{k}) \longrightarrow G(\overline{k}) \longrightarrow 1\]
of $\Gal(K/k)$-groups (where $K/k$ is the smallest finite field extension which splits both $G$ and $\mmu_{N^2}$) obtained from Lemma \ref{algoB}, the sublist of central extensions which represent the images of elements of $\Br^0_{nr}(\SL_{n,k}/G)$ in $\mathrm{Ext}^c_k(G(\overline{k}),\Q/\Z(1))$ (via the isomorphism given by Proposition \ref{normaliseextension}). 
Let us fix an extension $(E)$ provided by Lemma \ref{algoB}. We may apply Theorem \ref{nonramificationextension}, so that we have to check whether:
\begin{enumerate}
    \item[(i)] The extension $(E)$ splits (as an extension of abstract groups) when pulled back to any bicyclic subgroup of $G(\overline{k})$;
    \item[(ii)] For any $\sigma\in\Gamma_k\backslash\{1\}$ and for any two elements $\tau,\gamma\in G(\overline{k})$ such that $\gamma\mkern2.5mu {}^\sigma \mkern-2.5mu\tau \gamma^{-1}=\tau^{\chi(\sigma)}$, there exists a lift of abstract groups $\psi:\langle \tau\rangle\to E(\overline{k})$ and an element $e_\gamma\in E(\overline{k})$ which maps to $\gamma$ and verifies $e_\gamma\mkern1mu {}^\sigma \mkern-1mu\psi(\tau) e_\gamma^{-1}=\psi(\tau)^{\chi(\sigma)}$.
\end{enumerate}
By construction, $(E)$ verifies condition (i). Moreover since its class in $\mathrm{Ext}^c_k(G(\overline{k}),\Q/\Z(1))$ lifts to an element of $\mathrm{Ext}^c_k(G(\overline{k}),\mmu_N)$ for which the $\Gamma_k$-action factors through $\Gal(K/k)$, then condition (ii) only needs to be checked for $\sigma\in\Gal(K/k)$, that is, we fix any such $\sigma$ and any cyclic subgroup $C\subset G(\overline{k})$,  we list all of the possible (abstract) lifts $C\to E(\overline{k})$ and check which ones satisfy (ii). Indeed, the action of $\sigma\in\Gamma_k$ on $e_\gamma\psi(\tau)e_\gamma^{-1}$ here only depends on the image of $\sigma$ in $\Gal(K/k)$ since, by the very definition of the amalgamated sum, the element $e_\gamma\psi(\tau)e_\gamma^{-1}$ is represented by the sum of two elements of $\widetilde{E}(\overline{k})$ and $\mmu_{N^2}$ respectively. We can therefore decide whether or not $(E)$ represents an unramified class in a finite number of steps. Varying the choice of $(E)$ and applying the same sequence of arguments returns a finite list of central extensions representing all of the unramified classes. 

To check which elements of the list represent the same unramified Brauer class, it amounts to verifying which extensions provided by Lemma \ref{algoB} split when pushing forward along $\mmu_N\hookrightarrow\Q/\Z(1)$. This is once again achieved in a finite number of steps: the image of an element of order dividing $N$ by a morphism $G(\overline{k})\to \Q/\Z(1)\star_{\mmu_N}\widetilde{E}(\overline{k})$ is represented by the sum of two elements of $\widetilde{E}(\overline{k})$ and $\mmu_{N^2}$ respectively. 

Finally, in order to get the group law on $\Brnr^0(\SL_{n,k}/G)$ it suffices to understand the description of the sum of two central extensions of $\Gamma_k$-groups, which is classical and for instance described in \cite[Exposés III \& XVII]{SGA3}. 
\end{proof}

\section{Computing Brauer-Manin sets}
Let $G$ be a finite $k$-group together with an embedding of $k$-groups $G\hookrightarrow\SL_{n,k}$ for some suitable $n\ge 1$. We consider the problem of computing the Brauer-Manin obstruction to weak approximation on $\SL_{n,k}/G$ when $k$ is a number field. By \cite[Corollary 1.5]{harpazwittenbergsupersolvable}, this totally determines the Grunwald problem for $G$ as soon as the latter is supersolvable (e.g.\ constant and nilpotent).
\subsection{An algorithm to compute the Brauer-Manin obstruction} Thanks to \cite[Theorem 6.1, Theorem 6.2]{lucchiniunramifiedbrauer}, we know that the places that could produce a Brauer-Manin obstruction are finitely many, namely the ones that divide the order of $G(\overline{k})$ or ramify in the minimal field extension which splits $G$. We can therefore provide an algorithm to compute the Brauer-Manin set:
\begin{thm}\label{algorithmebrauermanin}
 There exists an algorithm that, given the datum of a finite $k$-group $G$ and a number field $k$, computes the Brauer-Manin obstruction to weak approximation on $\SL_{n,k}/G$ for any embedding $G\hookrightarrow\SL_{n,k}$ for some $n\ge1$.
\end{thm}
\begin{rmk}
    Let us make the statement of Theorem \ref{algorithmebrauermanin} more precise. The algorithm described in the proof of the latter produces a finite set of places $S\subset \Omega_k$ (which we can explicitely describe thanks to \cite[Theorem 6.1, Theorem 6.2]{lucchiniunramifiedbrauer}) and a subset $E_{BM}$ of the finite set $\prod_{v\in S}\rmH^1(k_v,G)$ that exactly determines the Brauer-Manin set of $\SL_{n,k}/G$ within $\prod_{v\in\Omega_k}(\SL_{n,k}/G)(k_v)$, in the sense that a collection $(P_v)_v$ of local points in the latter product lies in the Brauer-Manin set if and only if its image in $\prod_{v\in S}\rmH^1(k_v,G)$ lands into the finite subset $E_{BM}$.
\end{rmk}
\begin{proof} Start by computing the unramified Brauer group of $\SL_{n,k}/G$ by following the algorithm described in the proof of Theorem \ref{algorithmenonramifie}. This is achieved in a finite number of steps, and it returns a finite list $B\subset \mathrm{Ext}^c_k(G(\overline{k}),\Q/\Z(1))$ of (classes of) central extensions of $\Gamma_k$-groups which represent the unramified Brauer classes. 

Write a list $S$ of the finitely many places $v\in \Omega_k$ that divide $|G(\overline{k})|$ or ramify in the minimal field extension $E/k$ which splits $G$. By Lemma \ref{pointscalculables}, we know that the pointed sets $(\SL_{n,k}/G)(k_v)$ for $v\in S$ are finite and effectively computable. Thanks to Lemma \ref{comparaison}, for a given place $v\in  S$ the computation of the local Brauer-Manin pairing
\[\begin{tikzcd}
	{\Brnr^0(\SL_{n,k}/G)\times \rmH^1(k_v,G)} && {} & {\Br(k_v)}
	\arrow["{(\alpha,P_v)\mapsto \mathrm{ev}_\alpha(P_v)}", from=1-1, to=1-4]
\end{tikzcd}\]
can be conducted as follows: if a Brauer class $\alpha\in \Brnr^0(\SL_{n,k}/G)$ is represented by a central extension of $\Gamma_k$-groups
\[(E)~~:~~1\longrightarrow \Q/\Z(1)\longrightarrow E(\overline{k})\longrightarrow G(\overline{k})\longrightarrow 1,\]
then we can consider the boundary map $\alpha^*:\rmH^1(k_v,G)\to \rmH^2(k_v,\Q/\Z(1))= \Br(k_v)$ in the long exact Galois cohomology sequence associated to $(E)$. For every $P_v\in \rmH^1(k_v,G)$, we then have $\mathrm{ev}_\alpha(P_v)=\alpha^*(P_v)\in \Br(k_v)$.

The algorithm hence reads as follows. If $v\in \Omega_k\setminus S$, then \cite[Theorem 6.1,~Theorem~6.2]{lucchiniunramifiedbrauer} shows that the local evaluation pairing $ \Brnr^0(\SL_{n,k}/G)\times(\SL_{n,k}/G)(k_v)\to\Br(k_v)$ is constant and thus trivial, so there is nothing to do here. It therefore remains to compute local evaluations at the (potentially) bad places. Write $B=\{\alpha_1,\hdots,\alpha_r\}$, $S=\{v_1,\hdots,v_n\}$ and $\rmH^1(K_{v_i},G)=\{P_{v_i}^1,\hdots,P_{v_i}^{s_i}\}$ for $i\in \{1,\hdots,n\}$. Fix $\alpha_1$, then fix $v_1$, and compute $\mathrm{ev}_\alpha(P_{v_1}^j)$ for $j\in\{1,\hdots,{s_1}\}$. Once done, move on to $v_2$ and compute $\mathrm{ev}_\alpha(P_{v_2}^j)$ for $j\in\{1,\hdots,s_2\}$, then move on to $v_3$, \emph{etc}. Once that all of the $v_i$'s are covered, replace $\alpha_1$ by $\alpha_2$ and apply the same process, then replace $\alpha_2$ by $\alpha_3$, \emph{etc}. The algorithm eventually stops since all of the sets involved are finite.
\end{proof}

In the proof, we used the following:

\begin{lem}\label{pointscalculables}
Let $v\in \Omega_k$ be an arbitrary place. The pointed set $(\SL_{n,k}/G)(k_v)/\SL_{n}(k_v)$ is finite and effectively computable.\end{lem}
\begin{proof} For any field extension $L/k$, the elements of the pointed set $(\SL_{n,k}/G)(k_v)/\SL_{n}(k_v)$ correspond to elements of the pointed set $\rmH^1(L,G)$ thanks to Theorem~\ref{palschlank}. We hence need to compute the groups $\rmH^1(k_v,G)$ for all the finitely many places $v\in \Omega_k$ that divide $|G(\overline{k})|$ or may ramify in the minimal field extension $E/k$ which splits $G$. Let us start with the case where $G$ is constant. We then have $\rmH^1(k_v,G)={\Hom(\Gamma_{k_v},G)/(\text{conjugacy}})$; these groups are finite, since they classify the isomorphism classes of Galois extensions of $k_v$ with group $H\subset G$. The structure of local Galois groups \cite[Theorem 7.5.14]{neukirchschmidtwingberg} forces these isomorphism classes to be finitely many. If $G$ is now assumed to be non-necessarily constant, the set $\rmH^1(k_v,G)$ classifies $G$-torsors over $\Spec(k_v)$. If we let $E_v=E k_v$, then $G$ becomes constant over $E_v$ and a $G$-torsor over $\Spec(k_v)$ corresponds (up to isomorphism) to a $G$-torsor over $\Spec(E_v)$ which carries a descent datum. There are finitely many possible descent data since a $G$-torsor over $\Spec(E_v)$ admits only finitely many automorphisms. This readily implies the finiteness of $\rmH^1(k_v,G)$ in the general case, as well as its computability.
\end{proof}
\section{Examples}
\rm
We now describe some applications of the algorithm discussed in the proof of Theorem~\ref{algorithmenonramifie} in practical cases. For the sake of simplicity we only consider constant groups. Fix a field $k$ of characteristic $0$ and a constant $k$-group $G$ embedded into $\SL_{n,k}$ for some $n\ge 1$. 

\subsection{An example of transcendental Brauer-Manin obstruction}\label{paragraphDLAN} The first instance of Brauer-Manin obstruction to weak approximation for homogeneous spaces of linear groups given by a transcendental Brauer class was exhibited by Demarche, Lucchini-Arteche and Neftin in \cite[Example 5.3]{demarchelucchinineftin} (actually this is still the only known example together with the one constructed by Demeio in \cite{Dem22}). In \emph{loc.\ cit.}\, the authors derive a defect of weak approximation on $\SL_{n,k}/G$ for some semidirect product $G$ of abelian $p$-groups and some field $K$ which contains all the roots of unity of order dividing $\mathrm{exp}(G)$, which has to be explained by a Brauer-Manin obstruction thanks to Harari's theorem, cf.\ \cite[Théorème 1]{harariquelques}. Moreover this obstruction has to come from the only place of $k$ dividing $p$ since $G$ has $p$-power order \cite[Theorem 6.1, Theorem 6.2]{lucchiniunramifiedbrauer} and it has to be transcendental since the normalised algebraic subgroup of $\Brnr(\SL_{n,k}/G)$ vanishes in this case. However, their proof does not make the local evaluation providing a failure of weak approximation explicit, let alone the unramified Brauer class associated to it. We aim to exhibit a transcendental Brauer class which obstructs weak approximation on the homogeneous space considered in \cite[Example 5.3]{demarchelucchinineftin}, thereby removing the crucial dependence of their example on Harari's theorem.

To do this, we will begin by determining the unramified Brauer group of $\SL_{n,k}/G$ for every group $G$ which is a semidirect product of abelian groups, under the assumption that $k$ contains all the roots of unity of order $\mathrm{exp}(N)$.

\subsubsection{Bogomolov multipliers of semidirect products of abelian groups} Let $G=N\rtimes Q$ be the semidirect product of an abelian group $Q$ by another abelian group $N$. Recall that $B_0(G)$ can be alternatively described as the subgroup of $\rmH^2(G,\Q/\Z)$ consisting of classes which have trivial restrictions with respect to all abelian subgroups of $G$, see \cite[Lemma 3.10]{bogomolov}. We show that $B_0(G)$ can be identified with a subgroup of $\rmH^1(Q,\widehat{N})$:
\begin{prop}\label{bogomolovsemidirect}
Let $G=N\rtimes Q$ be a semidirect product of abelian groups. Then there is a natural isomorphism $\Sha^1_\mathrm{bic}(Q,\widehat{N})\simeq B_0(G)$.
\end{prop}
\begin{proof}
Consider the Hochschild-Serre spectral sequence 
\[E_2^{p,q}=\rmH^p(Q,\rmH^q(N,\Q/\Z))\Rightarrow \rmH^{p+q}(G,\Q/\Z).\]
As $G$ is a split extension of $Q$ by $N$, the prescribed section provides an exact sequence of lower terms:
\[0\longrightarrow \rmH^1(Q,\widehat{N})\longrightarrow \ker[\rmH^2(G,\Q/\Z)\longrightarrow \rmH^2(Q,\Q/\Z)]\longrightarrow \rmH^0(Q,\rmH^2(N,\Q/\Z)).\]
We thus see that $B_0(G)$ identifies to a subgroup of $\rmH^1(Q,\widehat{N})$. Let $B\subset Q$ be a bicyclic subgroup. We may pull back the split extension 
\[1\longrightarrow N\longrightarrow G\longrightarrow Q\longrightarrow 1\]
along the inclusion $B\subset Q$; the resulting subgroup $H\subset G$ is a semidirect product of $B$ by $N$ (the section being induced by the section of the above sequence). By the functoriality of the Hochschild-Serre spectral sequence, we hence obtain a commutative diagram with exact rows:
\[ 
    \begin{tikzcd} 
         0 \ar[r] & \rmH^1(Q,\widehat{N}) \ar[d] \ar[r] & \ker[\rmH^2(G,\Q/\Z)\to\rmH^2(Q,\Q/\Z)] \ar[d] \\
        0 \ar[r] & \rmH^1(B,\widehat{N}) \ar[r]& \ker[\rmH^2(H,\Q/\Z)\to\rmH^2(B,\Q/\Z)]
    \end{tikzcd}.
\]
If $\alpha\in\rmH^1(Q,\widehat{N})$ gives rise to a class in $B_0(G)\subset \ker[\rmH^2(G,\Q/\Z)\to\rmH^2(Q,\Q/\Z)]$, then it maps to $B_0(H)\subset \ker[\rmH^2(H,\Q/\Z)\to\rmH^2(B,\Q/\Z)]$ thanks to the functoriality of the Bogomolov multiplier. Now, by \cite[Lemma 3.2]{KreschTschinkel24}, the group $B_0(H)$ is trivial; since $\rmH^1(B,\widehat{N})$ injects into $\ker[\rmH^2(H,\Q/\Z)\to\rmH^2(B,\Q/\Z)]$ then this forces the image of $\alpha$ in $\rmH^1(B,\widehat{N})$ to be zero. This proves that $B_0(G)$ injects into $\Sha^1_\mathrm{bic}(Q,\widehat{N})$. 

Conversely, suppose that $\alpha$ lies in the latter group; we want to check whether its image in $\rmH^2(G,\Q/\Z)$ dies along the restriction to any bicyclic subgroup of $G$. Let $H\subset G$ be a bicyclic subgroup and let $B\subset Q$ be its image along the projection $H\hookrightarrow G\to Q$, which is bicyclic as well. Then $H$ is contained in the semidirect product $N\rtimes B$ (for the obvious action of $B$ on $N$). As before, the composite map $\rmH^1(Q,\widehat{N})\to \rmH^2(G,\Q/\Z)\to \rmH^2(N\rtimes B,\Q/\Z)$ factors through the injection $\rmH^1(B,\widehat{N})\hookrightarrow\rmH^2(N\rtimes B,\Q/\Z)$. Since $\alpha$ dies in $\rmH^1(B,\widehat{N})$ by assumption, then its image in $\rmH^2(G,\Q/\Z)$ vanishes in $\rmH^2(N\rtimes B,\Q/\Z)$ and hence also in $\rmH^2(H,\Q/\Z)$, as desired. 
\end{proof}

\subsubsection{Unramified Brauer groups for semidirect products of abelian groups} One can describe the elements of $\rmH^1(Q,\widehat{N})$ more precisely. Indeed, if $[a]\in\rmH^1(Q,\Hom(N,\Q/\Z))$ is a class represented by a $1$-cocycle $a:Q\to\Hom(N,\Q/\Z)$, then this cocycle determines a central extension (in a unique way)
\[1\longrightarrow \Q/\Z\longrightarrow (\Q/\Z\times N)\rtimes_a Q\longrightarrow N\rtimes Q\longrightarrow 1,\]
where the action of $Q$ on $\Q/\Z\times N$ is given by the formula
\[^q (\lambda,n):=(\lambda+a(q^{-1})(n),{^q n}).\]

Let $k$ be a field of characteristic $0$ and let $G\hookrightarrow \SL_{n,k}$ be an embedding of $k$-groups. As the natural map $\Br^0(\SL_{n,k}/G)\to \rmH^2(G,\Q/\Z)$ sends $\Brnr^0(\SL_{n,k}/G)$ to $B_0(G)$ then this shows, by combining with Proposition \ref{bogomolovsemidirect}, that any unramified Brauer class must be represented by a central extension of $\Gamma_k$-groups which, as an extension of abstract groups, is of the form described above. If we further assume that $k$ contains sufficiently many roots of unity, then we can precisely determine the entire unramified Brauer group of $\SL_{n,k}/G$:

\begin{prop}\label{actionsemidirect}
    Let $G=N\rtimes Q$ be an arbitrary semidirect product of abelian groups and let $k$ be a field of characteristic $0$ such that $\mmu_{\mathrm{exp}(N)}\subset k$. If $G\hookrightarrow\SL_{n,k}$ is any embedding of $k$-groups, then the natural map $\Brnr^0(\SL_{n,k}/G)\to B_0(G)\simeq \Sha^1_\mathrm{bic}(Q,\widehat{N})$ is onto. More precisely, any $1$-cocycle $a\in Z^1(Q,\widehat{N})$ such that $[a]\in \Sha^1_\mathrm{bic}(Q,\widehat{N})$ determines an unramified Brauer class $\alpha\in \Brnr^0(\SL_{n,k}/G)$ that is represented by the central extension of $\Gamma_k$-groups:
        \[ 
    \begin{tikzcd} 
        1 \ar[r] & \Q/\Z(1) \ar[r] & (\Q/\Z(1)\times N)\rtimes_{a} Q \ar[r] & N\rtimes Q \ar[r]& 1,
    \end{tikzcd}
\]
where $\Gamma_k$ acts on $(\Q/\Z(1)\times N)\rtimes_a Q$ only through the cyclotomic character on $\Q/\Z(1)$.
\end{prop}

We notably obtain the following special case of the above statement:
\begin{cor}
    Keeping the above notation, and assuming that $\Br_{\mathrm{nr},1}^0(\SL_{n,k}/G)=0$ (e.g.\ if $\mmu_{\exp(G)}\subset k$), then $\Brnr^0(\SL_{n,k}/G)\simeq \Sha^1_\mathrm{bic}(Q,\widehat{N})$. More precisely, any unramified Brauer class $\alpha\in \Brnr^0(\SL_{n,k}/G)$ is represented by a central extension of $\Gamma_k$-groups:
        \[ 
    \begin{tikzcd} 
        1 \ar[r] & \Q/\Z(1) \ar[r] & (\Q/\Z(1)\times N)\rtimes_{a} Q \ar[r] & N\rtimes Q \ar[r]& 1,
    \end{tikzcd}
\]
where $a\in Z^1(Q,\widehat{N})$ is a $1$-cocycle whose class determines an element $[a]\in\Sha^1_\mathrm{bic}(Q,\widehat{N})$ in a unique way, and $\Gamma_k$ acts on $(\Q/\Z(1)\times N)\rtimes_a Q$ only through the cyclotomic character on $\Q/\Z(1)$. 
\end{cor}

\begin{proof} Denote $e:=\mathrm{exp}(N)$. Thanks to Proposition \ref{bogomolovsemidirect}, we have that $B_0(G)\simeq \Sha^1_\mathrm{bic}(Q,\widehat{N}),$ and the pushforward map $\Sha^1_\mathrm{bic}(Q,\Hom(N,\Z/e))\to\Sha^1_\mathrm{bic}(Q,\widehat{N})$ associated to the inclusion $\Z/e\subset \Q/\Z$ is an isomorphism. As in the proof of Proposition \ref{nonramifieracinesunite}, we can fix an isomorphism of $\Gamma_k$-groups $\mmu_e\simeq \Z/e$ and consider the following commutative diagram:
\[
    \begin{tikzcd} 
          \mathrm{Ext}^c_k(G,\mmu_e) \ar[d, twoheadrightarrow] \ar[r] & \rmH^2(G,\Z/e)\ar[d, twoheadrightarrow]  \\
        \mathrm{Ext}^c_k(G,\Q/\Z(1))[e] \ar[r] & \rmH^2(G,\Q/\Z)[e]
    \end{tikzcd}
\]
where the top horizontal map admits an obvious section and $\Brnr^0(\SL_{n,k}/G)[e]$ injects into $\mathrm{Ext}_k^c(G,\Q/\Z(1))[e]$ (resp.\ $B_0(G)\subset \rmH^2(G,\Q/\Z)[e]$). Starting from a class $[a]\in \Sha^1_\mathrm{bic}(Q,\widehat{N})$, we may thus construct a class $\alpha\in \mathrm{Ext}_k^c(G,\Q/\Z(1))[e]$ that maps to $[a]$ and which is represented by a central extension that fits into a commutative diagram of extensions of $\Gamma_k$-groups:
\[ 
    \begin{tikzcd} 
         1 \ar[r] & \mmu_e \ar[d] \ar[r] & (\mmu_e\times N)\rtimes_{\widetilde{a}} Q \ar[d, "\varphi"]\ar[r] & G \ar[d, equal]\ar[r] & 1 \\
        1 \ar[r] & \Q/\Z(1) \ar[r]& (\Q/\Z(1)\times N)\rtimes_{a} Q \ar[r] & G \ar[r]& 1,
    \end{tikzcd}
\]
where the top row comes equipped with the trivial action of $\Gamma_k$ (here $\widetilde{a}\in Z^1(Q,\Hom(N,\mmu_e))$ is any representative of the inverse image of $[a]$ in $\Sha^1_\mathrm{bic}(Q,\Hom(N,\mmu_e))$). 

We may now apply Theorem \ref{nonramificationextension} to check that the class of the bottom row is unramified. We want to show that for any $\sigma\in \Gamma_k$ and any $\gamma,\tau\in G$ such that $\gamma\tau\gamma^{-1}=\tau^{\chi(\sigma)}$, there exists a lift of abstract groups $\psi:\langle \tau\rangle\to (\Q/\Z(1)\times N)\rtimes_a Q$ and an element $e_\gamma\in (\Q/\Z(1)\times N)\rtimes_a Q$ which maps to $\gamma$ and verifies $e_\gamma\mkern1mu {}^\sigma \mkern-1mu\psi(\tau)e_\gamma^{-1}=\psi(\tau)^{\chi(\sigma)}$. Note that it is sufficient to exhibit a section of abstract groups $\widetilde{\psi}:\langle\tau,\gamma\rangle\to (\Z/e\times N)\rtimes_a Q$: indeed, as $\varphi:(\mmu_e\times N)\rtimes_a Q\to (\Q/\Z(1)\times N)\rtimes_a Q$ is $\sigma$-equivariant, then letting $e_\gamma:=(\varphi\circ\widetilde{\psi})(\gamma)$ and $\psi:=\varphi\circ \widetilde{\psi}|_{\langle\tau\rangle}$ will provide the desired elements of $(\Q/\Z(1)\times N)\rtimes_a Q$, because
\[e_\gamma\mkern1mu {}^\sigma \mkern-1mu\psi(\tau)e_\gamma^{-1}=e_\gamma\psi(\tau)e_\gamma^{-1}=\varphi(\widetilde{\psi}(\gamma\tau\gamma^{-1}))=\phi(\widetilde{\psi}(\tau^{\chi(\sigma)}))=\psi(\tau)^{\chi(\sigma)}.\] As the class $[\widetilde{a}]$ of the top row of the above diagram lies in $\Sha^1_\mathrm{bic}(Q,\Hom(N,\Z/e))$, we may adapt the proof of Proposition \ref{bogomolovsemidirect}: if we denote by $H\subset G$ the subgroup generated by $\gamma$ and $\tau$, then its projection in $Q$ is a bicyclic subgroup $B$, and $H$ is thus contained in the semidirect product $N\rtimes B$ (for the obvious action of $B$ on $N$). On the other hand, the composite map $\rmH^1(Q,\Hom(N,\Z/e))\to \rmH^2(G,\Z/e)\to \rmH^2(N\rtimes B,\Z/e)$ factors through the injection $\rmH^1(B,\Hom(N,\Z/e))\hookrightarrow\rmH^2(N\rtimes B,\Z/e)$. Since $[\widetilde{a}]$ dies in $\rmH^1(B,\Hom(N,\Z/e))$ by assumption, then its image in $\rmH^2(G,\Z/e)$ vanishes in $\rmH^2(N\rtimes B,\Z/e)$ and hence also in $\rmH^2(H,\Z/e)$, as desired. Thus the image of $\alpha\in \mathrm{Ext}_k^c(G,\Q/\Z(1))$ in $\Br^0(\SL_{n,k}/G)$ is unramified.
\end{proof}

\subsubsection{Evaluation at local points}\label{evaluationpointslocaux} We keep the assumptions on $k$ and $G$. Suppose in addition that $k$ is a number field. Given the description of the unramified Brauer group of $\SL_{n,k}/G$ provided by Proposition \ref{actionsemidirect}, we will build upon the proof of \cite[Theorem 5.1]{demarchelucchinineftin} and make their result explicit and independent from Harari's work. 

Any transcendental Brauer class $\alpha\in \Brnr^0(\SL_{n,k}/G)$ provides a non-zero element of $B_0(G)\simeq \Sha^1_\mathrm{bic}(Q,\widehat{N})$; as Proposition \ref{palschlank} identifies the pointed sets $(\SL_{n,k_v}/G)(k_v)/\SL_{n}(k_v)$ and $\rmH^1(k_v,G)$, it is sufficient to construct a class $[c_v]\in\rmH^1(k_v,Q)$ whose fibre along the natural map $\rmH^1(k_v,G)\to \rmH^1(k_v,Q)$ contains a local element on which $\alpha$ evaluates non-trivially, hence providing a failure of weak approximation on $\SL_{n,k}/G$. 
\begin{prop}\label{obstructiontranscendante}
    Let $k$ be a number field such that $\mmu_{\mathrm{exp}(N)}\subset k$. For any transcendental unramified Brauer class $\alpha\in\Brnr^0(\SL_{n,k}/G)$ and for any place $v\in \Omega_k$ that divides the order of $G$, given an epimorphism $\Gamma_{k_v}\twoheadrightarrow Q$, there exists a local point $P_v\in (\SL_{n,k}/G)(k_v)$ such that $\mathrm{ev}_\alpha(P_v)\neq 0$. Consequently, the class $\alpha$ produces a transcendental Brauer-Manin obstruction to weak approximation on $\SL_{n,k}/G$.
\end{prop}
\begin{proof} Suppose that such a place $v\in\Omega_k$ exists (note that if $Q$ is a $p$-group, the structure of the multiplicative group $k_v^*$ \cite[Chapter II, Proposition 5.7]{NeukirchANT} forces $Q$ to be a Galois group over $k_v$ as soon as $k$ is large enough) and denote $c_v:\Gamma_{k_v}\twoheadrightarrow Q$ the prescribed epimorphism. Twisting the exact sequence underlying the semidirect product structure of $G$ by $c_v\in\Hom(\Gamma_{k_v},Q)\simeq\rmH^1(k_v,Q)$ yields a split short exact sequence of $\Gamma_{k_v}$-groups
\[ 
1\longrightarrow \mkern2mu {}_{c_v} \mkern-2mu N\longrightarrow \mkern2mu {}_{c_v} \mkern-2mu G\longrightarrow \mkern2mu {}_{c_v} \mkern-2mu Q\longrightarrow 1
\]
(where ${_{c_v} Q}=Q$ and $\Gamma_{k_v}$ acts on $\mkern2mu {}_{c_v} \mkern-2mu N$ as ${^\sigma n}:=c_v(\alpha)\cdot n\cdot c_v(\alpha)^{-1}\in \mkern2mu {}_{c_v} \mkern-2mu N$) which provides an injection of pointed sets  $\rmH^1(k_v,\mkern2mu {}_{c_v} \mkern-2mu N)\hookrightarrow\rmH^1(k_v,\mkern2mu {}_{c_v} \mkern-2mu G)$ after taking the associated long exact cohomology sequence. Moreover there exists a twisting bijection $\rmH^1(k_v,\mkern2mu {}_{c_v} \mkern-2mu G)\to \rmH^1(k_v,G)$, see \cite[ I.5.3, Proposition 35bis]{serrecg}. We thus obtain a composite map $\theta:\rmH^1(k_v,\mkern2mu {}_{c_v} \mkern-2mu N)\to \rmH^1(k_v,G)$. On the other hand, as $c_v$ was chosen to be surjective, then by inflation it induces an injective morphism $c_v^*:\rmH^1(Q,\widehat{N})\to \rmH^1(k_v,\mkern2mu {}_{c_v} \mkern-2mu \widehat{N})$. More precisely, this morphism takes an element $[a]\in\rmH^1(Q,\widehat{N})\hookrightarrow\rmH^2(G,\Q/\Z)$ represented by a central extension 
\[1\longrightarrow \Q/\Z(1)\longrightarrow (\Q/\Z(1)\times N)\rtimes_a Q\longrightarrow N\rtimes Q\longrightarrow 1\]
and sends it to the class associated to the pullback of this extension along the composition $\mkern2mu {}_{c_v} \mkern-2mu N\hookrightarrow \mkern2mu {}_{c_v} \mkern-2mu N\rtimes\Gamma_{k_v}\twoheadrightarrow N\rtimes Q$ (where the second map is given by $(n,\sigma)\mapsto(n,c_v(\sigma))$), equipped with the $\Gamma_{k_v}$-action given by ${^\sigma(\lambda,n)}:=(\lambda^{\chi(\sigma)}+a(c_v(\sigma)^{-1})(n),{^{c_v(\sigma)}n})$ (this is well defined since $\Gamma_{k_v}$ acts on $_{c_v}\widehat{N}=\Hom(_{c_v}N,\mmu_{\mathrm{exp(N)}})$ only through the prescribed action of $Q$ on $N$ via the cocycle $c_v$, as it acts trivially on $\mmu_{\mathrm{exp(N)}}$).

We may thus consider a non-zero class $[a]\in \rmH^1(Q,\widehat{N})$ and let $\alpha\in\Brnr^0(\SL_{n,k}/G)$ be any element whose image along the natural map $\Brnr^0(\SL_{n,k}/G)\to B_0(G)\simeq \Sha^1_\mathrm{bic}(Q,\widehat{N})$ is $[a]$ (such an element exists thanks to Proposition \ref{actionsemidirect}).

As $c_v^*([a])\neq0\in\rmH^1(k_v,\mkern2mu {}_{c_v} \mkern-2mu \widehat{N})$, we can then apply local duality \cite[Theorem~7.2.6]{neukirchschmidtwingberg}; this implies the existence of a local point $P_{v}\in\rmH^1(k_{v},\mkern2mu {}_{c_v} \mkern-2mu N)$ on which $c_v^*([a])$ evaluates non-trivially. Lemma \ref{commutatif} then shows that $P_v$ yields a local point $P_v'\in \rmH^1(k_v,G)$ on which $\alpha$ evaluates non-trivially as well, hence providing the desired transcendental Brauer-Manin obstruction to weak approximation. 
\end{proof}

In the proof, we used the following:
\begin{lem}\label{commutatif}
 There exists an adjunction diagram: 
\[\begin{tikzcd}
	{\Brnr^0(\SL_{n,k}/G)} & \times & {\rmH^1(k_v,G)} & {\Br(k_v)} \\
	{\rmH^1(k_v,\mkern2mu {}_{c_v} \mkern-2mu \widehat{N})} & \times & {\rmH^1(k_v,\mkern2mu {}_{c_v} \mkern-2mu N)} & {\Br(k_v)},
	\arrow["\iota_{c_v}^*", from=1-1, to=2-1]
	\arrow[from=1-3, to=1-4]
	\arrow[equals, from=1-4, to=2-4]
	\arrow["\theta", hook, from=2-3, to=1-3]
	\arrow[from=2-3, to=2-4]
\end{tikzcd}\]
where the leftmost vertical map is given by the composition of the natural surjective map $\Brnr^0(\SL_{n,k}/G)\twoheadrightarrow \Sha^1_\mathrm{bic}(Q,\widehat{N})$ from Proposition \ref{actionsemidirect} and $c_v^*:\rmH^1(Q,\widehat{N})\to \rmH^1(k_v,\mkern2mu {}_{c_v} \mkern-2mu \widehat{N})$, the top horizontal arrow is the Brauer-Manin pairing at the place $v$, and the bottom horizontal one is the local duality pairing. Moreover, $\iota_{c_v}^*$ is trivial on $\Br_{\mathrm{nr},1}^0(\SL_{n,k}/G)$ and induces an injection $\overline{\iota_{c_v}^*}:\Brnr^0(\SL_{n,k}/G)/\Br_{\mathrm{nr},1}^0(\SL_{n,k}/G)\hookrightarrow\rmH^1(k_v,\mkern2mu {}_{c_v} \mkern-2mu \widehat{N})$.
\end{lem}
\begin{proof}
Start with a class $\alpha\in \Brnr^0(\SL_{n,k}/G)$ and view it as an element of $\Br(\SL_{n,k}/G)\simeq \rmH^2(G\times\Gamma_k,\Q/\Z(1))$. There is a morphism of semidirect products $\iota_{c_v}:\mkern2mu {}_{c_v} \mkern-2mu N\rtimes \Gamma_{k_v}\to G\times \Gamma_k$ given by $(\mathrm{Id}:\mkern2mu {}_{c_v} \mkern-2mu N\to N)\times(c_v\times\mathrm{Id}:\Gamma_{k_v}\to Q\times\Gamma_k)$, that induces a morphism $\iota_{c_v}:\mkern2mu {}_{c_v} \mkern-2mu N\to G$ which is $c_v$-twisted-equivariant for the actions of $\Gamma_{k_v}$ on $\mkern2mu {}_{c_v} \mkern-2mu N$ and $\Gamma_k$ on $G$ respectively, in the sense that $\iota_{c_v}({^\sigma n})=c_v(\sigma)\cdot \iota_{c_v}(n)\cdot c_v(\sigma)^{-1}\in N\subset G$ for any $n\in \mkern2mu {}_{c_v} \mkern-2mu N$ and any $\sigma\in \Gamma_{k_v}$. Pulling back along this morphism thus gives a map $\iota_{c_v}^*:\rmH^2(G\times\Gamma_k,\Q/\Z(1))\to \rmH^2(\mkern2mu {}_{c_v} \mkern-2mu N\rtimes\Gamma_{k_v},\Q/\Z(1))$, and we may compose it with the split surjection $\rmH^2(\mkern2mu {}_{c_v} \mkern-2mu N\rtimes\Gamma_{k_v},\Q/\Z(1))\to \mathrm{Ext}_{k_v}^c(\mkern2mu {}_{c_v} \mkern-2mu N,\Q/\Z(1))$ from Proposition \ref{normaliseextension}. Actually, since $\alpha$ was chosen to be unramified, then its image along $\Brnr^0(\SL_{n,k}/G)\to \mathrm{Ext}_{k_v}^c(\mkern2mu {}_{c_v} \mkern-2mu N,\Q/\Z(1))$ is split as an extension of abstract groups, hence uniquely determines an element of $\rmH^1(k_v,\mkern2mu {}_{c_v} \mkern-2mu\widehat{N})$: this is clear by looking at the corresponding diagram of pullback extensions:
\[\begin{tikzcd}
	1 & {\Q/\Z(1)} & {(\Q/\Z(1)\times N)\rtimes Q} & {N\rtimes Q} & 1 \\
	1 & {\Q/\Z(1)} & {\Q/\Z(1)\times {\mkern2mu {}_{c_v} \mkern-2mu N}} & {{\mkern2mu {}_{c_v} \mkern-2mu N}} & 1
	\arrow[from=1-1, to=1-2]
	\arrow[from=1-2, to=1-3]
	\arrow[from=1-3, to=1-4]
	\arrow[from=1-4, to=1-5]
	\arrow[from=2-1, to=2-2]
	\arrow[equals, from=2-2, to=1-2]
	\arrow[from=2-2, to=2-3]
	\arrow[from=2-3, to=1-3]
	\arrow[from=2-3, to=2-4]
	\arrow["\iota_{c_v}"', from=2-4, to=1-4]
	\arrow[from=2-4, to=2-5]
\end{tikzcd}\]
where we require for all the vertical maps to be $c_v$-twisted-equivariant (here the $\Gamma_{k}$-action on $\big((\Q/\Z(1)\times N)\rtimes Q$ is implicitely prescribed, we write the underlying abstract groups in the above diagram for the sake of clarity). We hence obtain a map $\iota_{c_v}^*:\Brnr^0(\SL_{n,k}/G)\to \rmH^1(k_v,\mkern2mu {}_{c_v} \mkern-2mu \widehat{N})$. Again, the $\Gamma_{k_v}$-action on $\Q/\Z(1)\times \mkern2mu {}_{c_v} \mkern-2mu N$ in the middle row corresponds to a $\Gamma_{k_v}$-action on $\mkern2mu {}_{c_v} \mkern-2mu \widehat{N}=\Hom(\mkern2mu {}_{c_v} \mkern-2mu N,\mmu_{\mathrm{exp}(N)})$ which only depends, via the cocycle $c_v$, on the prescribed action of $Q$ on $N$ (as $\Gamma_{k_v}$ acts trivially on $\mmu_{\mathrm{exp}(N)}$). This shows that the $\Gamma_{k_v}$-action on the bottom row of the diagram only depends of the underlying class of abstract groups of the top extension by forgetting the $\Gamma_k$-action (i.e.\ seen as an element of $\rmH^1(Q,\widehat{N})\hookrightarrow\rmH^2(G,\Q/\Z)$). Since this class is precisely the image of $\alpha$ in $\rmH^1(Q,\widehat{N})$ via the map $\Brnr^0(\SL_{n,k}/G)\to\Sha^1_{\mathrm{bic}}(Q,\widehat{N})$, then this readily implies that the map $\iota:\Brnr^0(\SL_{n,k}/G)\to \rmH^1(k_v,\mkern2mu {}_{c_v} \mkern-2mu \widehat{N})$ coincides with 
\[\Brnr^0(\SL_{n,k}/G)\overset{}{\longrightarrow} \Sha^1_\mathrm{bic}(Q,\widehat{N})\overset{c_v^*}{\longrightarrow}\rmH^1(k_v,\mkern2mu {}_{c_v} \mkern-2mu\widehat{N}),\]
and hence is zero on $\Br_{\mathrm{nr},1}^0(\SL_{n,k}/G)$. As $\Brnr^0(\SL_{n,k}/G)/\Br_{\mathrm{nr},1}^0(\SL_{n,k}/G)\simeq \Sha^1_\mathrm{bic}(Q,\widehat{N})$ and $c_v^*$ is injective, then so is $\overline{\iota_{c_v}^*}:\Brnr^0(\SL_{n,k}/G)/\Br_{\mathrm{nr},1}^0(\SL_{n,k}/G)\to \rmH^1(k_v,\mkern2mu {}_{c_v} \mkern-2mu\widehat{N})$. It therefore remains to check that the following diagram commutes:
\[ 
    \begin{tikzcd} 
        \rmH^1(k_{v},\mkern2mu {}_{c_v} \mkern-2mu\widehat{N}) \ar[d, hookrightarrow, "\theta"'] \ar[r, "\mathrm{ev}_{\iota(\alpha)}"] & \Br(k_{v}) \ar[d, equal] \\
       \rmH^1(k_{v},{G}) \ar[r, "\mathrm{ev}_\alpha"] & \Br(k_{v}).
    \end{tikzcd}
\]
But this follows, thanks to the very construction of $\iota_{c_v}^*$, from Lemma \ref{comparaison} and the compatibility of local duality with twisting bijections, see Steps 4 and 5 in the proof of \cite[Theorem~5.1]{demarchelucchinineftin}.
\end{proof}

\begin{rmk}\label{remarqueDLAN} The proof of \cite[Theorem 5.1]{demarchelucchinineftin} relies on Poitou-Tate duality in order to derive a defect of weak approximation on $\SL_{n,k}/G$, as the construction of local classes that are globally unattainable in \emph{loc.\ cit.}\ cannot make use of an explicit unramified Brauer class in order to evaluate them. The slight improvement of our approach lies in the sole use of local duality at a fixed place, which is only made possible by the description of the unramified Brauer group obtained in Proposition \ref{actionsemidirect}. In particular, we derive a stronger statement: for any place $v\in\Omega_k$ that divides the order of $G$ and such that $Q$ can be realised as a Galois group over $k_v$, and for any transcendental unramified Brauer class $\alpha\in \Brnr^0(\SL_{n,k}/G)$, there exists a Brauer-Manin obstruction to weak approximation coming from a non-zero local evaluation of $\alpha$ at the place $v$.
\end{rmk}

\begin{rmk}
Suppose that $G=N\rtimes Q$ is an arbitrary semidirect product of abelian $p$-groups for which $B_0(G)\simeq\Sha^1_\mathrm{bic}(Q,\widehat{N})\neq 0$. Fix an embedding $G\hookrightarrow \SL_{n,k}$ for some number field $k$ over which every class of $B_0(G)$ descends, e.g.\ $k=\Q(\zeta_{\mathrm{exp(N)}})$. Thanks to Proposition \ref{obstructiontranscendante}, we see that any transcendental class $\alpha\in \Brnr^0(\SL_{n,k}/G)$ provides a Brauer-Manin obstruction to weak approximation as soon as $Q$ identifies with the Galois group of some local field extension $L^{(v)}/k_v$ for some place $v\in\Omega_k$ dividing $p$. This obstruction notably persists over any finite extension of $k$. Indeed, conditions (1), (3) and (4) in \emph{loc.\ cit.\ ibidem} are stable under taking finite extensions of $k$ (up to enlarging the set of places), and by local class field theory, as soon as the finite abelian $p$-group $Q$ is a Galois group with respect to some fixed place $v\in\Omega_k$, then the structure of the multiplicative group $k_v^*$ \cite[Chapter II, Proposition 5.7]{NeukirchANT} forces $Q$ to be a Galois group over any finite extension of $k_v$.
\end{rmk}

\subsubsection{An explicit example} Let $p$ be a prime number and consider the group $Q:=(\Z/p)^3$. Let $I$ be the augmentation ideal of $(\Z/p^3)[Q]$ and let $N:=\widehat{I}=\Hom(I,\Q/\Z)$. There is a natural action of $Q$ on $N$ and $G:=N\rtimes Q$ is a group of order $p^{3p^3}$. We let $K:=\Q(\zeta_{p^3})$. The group $I$ fits by definition into the following exact sequence of $Q$-modules:
\[0\longrightarrow I\longrightarrow \Z/p^3[Q]\longrightarrow \Z/p^3\longrightarrow0,\]
where $\Z/p^3[Q]\to\Z/p^3$ is the augmentation map. For any $G$-module $M$, let $\widehat{H}^0(G,M)$ denote the Tate cohomology group $\rmH^0(G,M)/N_G(M)$, where $N_G(M):=\{\sum_{g\in G} {^g m}~|~m\in M\}$. Since $\Z/p^3[Q]$ is an induced $Q$-module, then the associated long exact cohomology sequence and Shapiro's lemma provide an isomorphism $\rmH^1(Q,I)\simeq \widehat{\rmH}^0(Q,\Z/p^3)\simeq \Z/p^3$ (resp.\ $\widehat{\rmH}^0(B,\Z/p^3)\simeq \Z/|B|$ for any subgroup $B\subset Q$). In particular, we get: 
\[\Sha^1_\mathrm{bic}(Q,\widehat{N})\simeq\ker\bigg[\Z/|Q|\to \prod_{B\subset Q}\Z/|B|\bigg]=\Z/p,\]
where the product ranges over the set of bicyclic subgroups $B$ of $Q$. In particular, this group is generated by the class $p^2[a]\in \rmH^1(Q,\widehat{N})$, where $[a]$ is represented by the $1$-cocycle $a:Q\to\Hom(N,\Q/\Z),~ q\mapsto [q]-[1_Q]$. Applying Proposition \ref{actionsemidirect}, we thus obtain the following result:
\begin{cor}\label{descriptionbrauernonramifie}
    Any transcendental element of $\Brnr^0(\SL_{n,K}/G)$ is of the form $\alpha+\beta$ where $\beta\in\Br_{\mathrm{nr},1}^0(\SL_{n,K}/G)$ and $\alpha$ is the transcendental Brauer class of order $p$ represented by the extension of $\Gamma_K$-groups:

    \[ 
    \begin{tikzcd} 
        1 \ar[r] & \Q/\Z(1) \ar[r] & (\Q/\Z(1)\times N)\rtimes_{p^2a} Q \ar[r] & N\rtimes Q \ar[r]& 1,
    \end{tikzcd}
\]
where $\Gamma_K$ acts on $(\Q/\Z(1)\times N)\rtimes_{p^2a} Q$ only through the cyclotomic character on $\Q/\Z(1)$. 
\end{cor}
If now $v\in \Omega_K$ denotes a place that divides $p$, then $Q$ can be viewed as a quotient of $K_v^*$ and hence can be realised as a Galois group over $K_v$ thanks to local class field theory. Taking $c_v:\Gamma_{K_v}\twoheadrightarrow Q$ to be the corresponding morphism and applying Proposition~\ref{obstructiontranscendante} provides a local point $P_v\in (\SL_{n,K}/G)(K_v)$ on which $\alpha$ (and \emph{a fortiori} $\alpha+\beta$ for any $\beta\in \Br_{\mathrm{nr},1}^0(\SL_{n,K}/G)$) evaluates non-trivially. We thus derive a Brauer-Manin obstruction to weak approximation on $\SL_{n,K}/G$.
\begin{rmk}
    The transcendental Brauer-Manin obstructions constructed above cannot descend to $\SL_{n,\Q}/G$ if we assume that $p$ is odd. Indeed, the normalised unramified Brauer group $\Brnr^0(\SL_{n,\Q}/G)$ vanishes in this case \cite[Corollaire 5.7, (b)]{ctnonramifie}. More generally, the same remark applies to the transcendental Brauer-Manin obstructions constructed in §\ref{evaluationpointslocaux} when $G$ is an arbitrary semi-direct product of abelian $p$-groups with $p$ odd.
\end{rmk}
\subsection{Transcendental classes over the real numbers}\label{paragraphkunyavskii} Let $G$ be a finite group, let $k=\R$, and view $G$ as a constant $\R$-group. Fix an embedding $G\hookrightarrow\SL_{n,\R}$ for some $n\ge 1$. The first part of the following lemma was already known thanks to the work of Colliot-Thélène \cite{ctnonramifie}:

\begin{lem}\label{brauerréel}
    The group $\Brnr^0(\SL_{n,\R}/G)$ is killed by $2$. In particular if $G$ has odd order, or if its $2$-Sylow subgroups are abelian, then $\Brnr^0(\SL_{n,\R}/G)=0$. Moreover, the image of any class $\alpha\in \Brnr^0(\SL_{n,\R}/G)$ in $\mathrm{Ext}^c_\R(G,\Q/\Z(1))$ via the isomorphism of Proposition \ref{normaliseextension} lifts to an element of $\mathrm{Ext}_\R^c(G,\mmu_2)$.
\end{lem}
\begin{proof} Since $G$ is constant, we have an isomorphism $\Br(\SL_{n,\R}/G)\simeq \rmH^2(G,\R^*)$ (see e.g.\ the proof of \cite[Proposition 1.1]{demarchenilpotent}). As $\mmu_\infty(\R)=\mmu_2$, we may use the tautological exact sequence
\[1\longrightarrow \mmu_2\longrightarrow \R^*\longrightarrow \R^*/\mmu_2\longrightarrow 1\]
and the fact that $\R^*/\mmu_2$ is torsion-free and divisible to derive an isomorphism $\rmH^2(G,\R^*)\simeq \rmH^2(G,\mmu_2$ as in \cite[Théorème 5.5, (b)]{ctnonramifie}. In particular $\Brnr^0(\SL_{n,\R}/G)$ is killed by $2$ and we may lift any class in $\Brnr^0(\SL_{n,\R}/G)$ to the class of a central extension in $\mathrm{Ext}_\R^c(G,\mmu_2)$ using the Kummer sequence \cite[Exposé XVII, Proposition A.2.1, a)]{SGA3}. If $G$ has odd order, then $\rmH^2(G,\mmu_2=0$, so that $\Brnr^0(\SL_{n,\R}/G)=0$. Similarly, if $S$ is a $2$-Sylow subgroup of $G$, then by restriction-corestriction we have an injection $\Brnr^0(\SL_{n,\R}/G)\hookrightarrow \Brnr^0(\SL_{n,\R}/S)$. Thanks to e.g.\ \cite[Proposition 5.7]{lucchiniartechenonramifie}, the group $\Brnr^0(\SL_{n,\R}/S)$ is trivial as soon as $S$ is constant and abelian, hence the result. 
\end{proof}
\subsubsection{A shorter version of the algorithm of Theorem \ref{algorithmenonramifie}} Lemma \ref{brauerréel} allows us to shorten the algorithm of Theorem \ref{algorithmenonramifie} in the real case, and we may assume that $G$ is a $2$-group. The algorithm reads as follows:

\begin{enumerate}
    \item[(i)] To determine the transcendental classes in $\Brnr^0(\SL_{n,\R}/G)$ consider, for any non-zero class $\alpha\in B_0(G)$, a lift $\widetilde{\alpha}\in \rmH^2(G,\Z/2)$ and view it as an element of $\mathrm{Ext}_\R^c(G,\mmu_2)$ equipped with the trivial $\Gamma_\R$-action. Pushing the underlying extension forward along the inclusion $ \mmu_2\hookrightarrow\Q/\Z(1)$ naturally provides an extension of $\Gamma_\R$-groups which, as an extension of abstract groups, corresponds to $\alpha\in B_0(G)$. Let us denote by $\alpha'\in\mathrm{Ext}_\R^c(G,\Q/\Z(1))$ the class of this extension. To get all of the classes in $\Br^0(\SL_{n,\R}/G)\simeq \mathrm{Ext}_\R^c(G,\Q/\Z(1))$ which map to $\alpha\in B_0(G)$, we replace $\alpha'$ by $\alpha'+\beta$ for all $\beta\in\Br_1^0(\SL_{n,R}/G)\simeq\rmH^1(\R,\widehat{G}^\mathrm{ab})$.
    \item[(ii)] Apply Theorem \ref{nonramificationextension} to the obtained classes in $\Br^0(\SL_{n,\R}/G)\simeq \mathrm{Ext}_\R^c(G,\Q/\Z(1))$ to determine which ones are unramified. 
\end{enumerate}
\subsubsection{An
 example with a quasisimple group}\label{exempleréels} For instance, let us consider the group $\mathrm{PSL}_3(\mathbf{F}_4)$, which is perfect (i.e.\ satisfies $\mathrm{PSL}_3(\mathbf{F}_4)^\mathrm{ab}=0$); we can thus form its universal central extension \cite[Proposition 33.4]{aschbacher}:
 \[1\longrightarrow\rmH_2(\mathrm{PSL}_3(\mathbf{F}_4),\Z)\longrightarrow \widetilde{G}\longrightarrow \mathrm{PSL}_3(\mathbf{F}_4)\longrightarrow 1,\]
 where the universal coefficient theorem identifies the groups $\rmH^2(\mathrm{PSL}_3(\mathbf{F}_4),\Q/\Z)$ and $\Hom(\rmH_2(\mathrm{PSL}_3(\mathbf{F}_4),\Z),\Q/\Z)$ (in particular, the kernel of the above extension is canonically isomorphic to the Pontrjagin dual of $\rmH^2(\mathrm{PSL}_3(\mathbf{F}_4),\Q/\Z)$). By \cite[§5]{bogomolovsimple}, we have $\rmH^2(\mathrm{PSL}_3(\mathbf{F}_4),\Q/\Z)\simeq \Z/4\oplus \Z/12$. Pushing the above extension forward along the evaluation at an element $\gamma\in \rmH^2(\mathrm{PSL}_3(\mathbf{F}_4),\Q/\Z)$ of order $4$ yields a central extension 
\begin{equation}\label{extensionkunyavskiiintermediaire}1\longrightarrow \Z/4\longrightarrow G\longrightarrow \mathrm{PSL}_3(\mathbb{F}_4)\longrightarrow 1.\end{equation}
The middle group is therefore a $4$-cover of $\mathrm{PSL}_3(\mathbf{F}_4)$. By construction, we have an epimorphism $\widetilde{G}\twoheadrightarrow G$ whose kernel $Z$ coincides with the kernel of the evaluation map $\mathrm{ev}_\gamma:\Hom(\rmH^2(\mathrm{PSL}_3(\mathbf{F}_4),\Q/\Z),\Q/\Z)\to \Z/4$.

As it was remarked by Kunyavski\u{\i}, any such $4$-cover $G$ satisfies $B_0(G)=\Z/2$ \cite[Theorem 3.1]{kunyavskii}. More precisely, since $\widetilde{G}$ is a central extension of $G$ by $Z$, we may look at the exact sequence of lower terms in the Hochschild-Serre spectral sequence $E_2^{p,q}=\rmH^p(G,\rmH^q(Z,\Q/\Z))\Rightarrow\rmH^{p+q}(\widetilde{G},\Q/\Z)$:
\[\rmH^1(\widetilde{G},\Q/\Z)\longrightarrow \rmH^1(Z,\Q/\Z)\longrightarrow \rmH^2(G,\Q/\Z)\longrightarrow \rmH^2(\widetilde{G},\Q/\Z).\]
Since $\widetilde{G}$ is the universal central extension of a perfect group, then it is perfect and satisfies $\rmH^2(\widetilde{G},\Q/\Z)=0$, see \cite[Proposition 33.8]{aschbacher}. We thus see that $\rmH^2(G,\Q/\Z)\simeq \widehat{Z}\simeq \Z/12$. Pushing $\widetilde{G}$ out (viewed as an extension of $G$ by $Z$) along the inclusion $Z\hookrightarrow \Q/\Z$ gives a central extension 
\begin{equation}\label{extensionkunyavskii}
1\longrightarrow \Q/\Z\longrightarrow E\longrightarrow G\longrightarrow 1
\end{equation}
which represents a generator of $H^2(G,\Q/\Z)$ (the kernel of the pushforward map $\rmH^2(G,Z)\to\rmH^2(G,\Q/\Z)$ is a quotient of $\rmH^1(G,\Q/\Z)=\widehat{G}^\mathrm{ab}=0$).
\begin{lem}\label{braueralgebriquereel}
    We have $\Br_{1}^0(\SL_{n,\R}/G)=0$; as a consequence, we have isomorphisms $\Br^0(\SL_{n,\R}/G)\simeq \rmH^2(G,\Q/\Z)[2]\simeq\rmH^2(G,\Z/2)\simeq \Z/2$, and the generator of $\Br^0(\SL_{n,\R}/G)$ is represented by a central extension of $\Gamma_\R$-groups
    \begin{equation}\label{extensionkunyavskiigalois}1\longrightarrow \Q/\Z(1)\longrightarrow E(\C)\longrightarrow G\longrightarrow 1\end{equation}
    such that $E(\C)$ is equipped with a non-trivial $\Gamma_\R$-action and is isomorphic, as an abstract group, with the group $E$ in \eqref{extensionkunyavskii}.
\end{lem}
\begin{proof} As $\Br_1^0(\SL_{n,\R}/G)\simeq\rmH^1(\R,\widehat{G}^\mathrm{ab})$ and the group $G^\mathrm{ab}$ is trivial, we get the desired vanishing. In particular, $\Br^0(\SL_{n,\R}/G)$ injects into $\rmH^2(G,\Q/\Z)$, and by Lemma \ref{brauerréel} it identifies \emph{a fortiori} with a subgroup of $\rmH^2(G,\Q/\Z)[2]$. On the other hand, we may identify $\Br^0(\SL_{n,\R}/G)\simeq\mathrm{Ext}_\R^c(G,\Q/\Z(1))$ and consider the following commutative diagram as in the proof of Proposition \ref{nonramifieracinesunite}:
\[
    \begin{tikzcd} 
          \mathrm{Ext}^c_\R(G,\mmu_2) \ar[d, twoheadrightarrow] \ar[r] & \rmH^2(G,\Z/2)\ar[d, twoheadrightarrow]  \\
        \mathrm{Ext}^c_\R(G,\Q/\Z(1)) \ar[r] & \rmH^2(G,\Q/\Z)[2]
    \end{tikzcd}
\]
(here we use the fact that $\Z/2\simeq \mmu_2$ as $\Gamma_\R$-modules). The kernels of all the maps involved in this diagram are trivial since $G^\mathrm{ab}=0$. By looking at the evident section of $\mathrm{Ext}^c_\R(G,\mmu_2)\to \rmH^2(G,\Z/2)$, we see that the latter map is surjective, hence an isomorphism. Therefore, all the maps in the diagram are isomorphisms. This implies that the generator of $\Br^0(\SL_{n,\R}/G)\simeq\mathrm{Ext}_\R^c(G,\Q/\Z(1))$ is represented by an extension
\[1\longrightarrow \Q/\Z(1)\longrightarrow E(\C)\longrightarrow G\longrightarrow 1\]
which can be realised as the pushforward of the only constant non-trivial central extension of $G$ by $\Z/2$.
\end{proof}
Therefore, the only possibility for a class in $\Br^0(\SL_{n,\R}/G)$ to be unramified is that of the one represented by the central extension of $\Gamma_\R$-groups \eqref{extensionkunyavskiigalois}. 
\begin{prop}\label{kunyavskiitrivial}
   The class of $\Br^0(\SL_{n,\R}/G)$ represented by \eqref{extensionkunyavskiigalois} is ramified. Consequently, we have $\Brnr^0(\SL_{n,\R}/G)=0$. 
\end{prop}
\begin{proof}
    Let $\widetilde{\alpha}\in \rmH^2(G,\Z/2)\simeq \mathrm{Ext}_\R^c(G,\mmu_2)\simeq \Z/2$ be the generator, and suppose that its image $\alpha\in\mathrm{Ext}_\R^c(G,\Q/\Z(1))\simeq\Br^0(\SL_{n,\R}/G)$ is unramified. Then the image of the latter along the restriction
    \[\mathrm{Ext}_\R^c(G,\Q/\Z(1))\longrightarrow \mathrm{Ext}_\R^c(A,\Q/\Z(1))\]
    with respect to any abelian subgroup $A\subset G$ is trivial, thanks to \cite[Proposition~5.7]{lucchiniartechenonramifie}. However, the image of $\rmH^2(G,\Z/2)\simeq\mathrm{Ext}_\R^c(G,\mmu_2)$ in $\mathrm{Ext}_\R^c(A,\mmu_2)$ factors by construction through $\rmH^2(A,\Z/2)$, which only consists of classes of extensions of $\Gamma_\R$-groups equipped with the trivial action. As $\R$ does not contain $\mmu_4$, the image of any non-zero element of $\Hom(A,\mmu_2)$ in $\mathrm{Ext}_\R^c(A,\mmu_2)$ along the Kummer map is on the other hand represented by a central extension endowed with a non-trivial $\Gamma_\R$-action. More precisely, any element arising this way is represented by a central extension (E) that fits as the top row of a commutative diagram of extensions of $\Gamma_\R$-groups:
    \[ 
    \begin{tikzcd} 
      1 \ar[r] & \mmu_2 \ar[d] \ar[r] & E(\C) \ar[d]\ar[r] & A \ar[d]\ar[r] & 1 \\ 
        1 \ar[r] & \mmu_2 \ar[r]& \mmu_4 \ar[r] & \mmu_2 \ar[r]& 1,
    \end{tikzcd}
\]
where all vertical maps are surjective; thus, the non-trivial $\Gamma_\R$-action on $\mmu_4$ induces a non-trivial one on $E(\C)$. Since the image of the Kummer map is the kernel of the pushforward $\mathrm{Ext}_\R^c(A,\mmu_2)\to\mathrm{Ext}_\R^c(A,\Q/\Z(1))$, this shows that $\widetilde{\alpha}\in\rmH^2(G,\Z/2)\simeq \mathrm{Ext}_\R^c(G,\mmu_2)$ must die in $\rmH^2(A,\Z/2)\hookrightarrow\mathrm{Ext}_\R^c(A,\mmu_2)$ for any abelian subgroup $A$ of $G$, and hence must lie in $\Sha^2_\mathrm{ab}(G,\Z/2)$. For the group $G$ that we consider in §\ref{exempleréels}, a computation in GAP shows that $\Sha^2_\mathrm{ab}(G,\Z/2)=0$. This ensures that the class represented by \eqref{extensionkunyavskiigalois} cannot be unramified.
\end{proof}
\begin{rmk}
    If we let $k$ be a number field equipped with a real place (e.g.\ $k=\Q$), then any unramified Brauer class $\alpha\in\Brnr^0(\SL_{n,k}/G)$ is algebraic: indeed, the natural map $\Brnr^0(\SL_{n,k}/G)\to B_0(G)$ factors through $\Brnr^0(\SL_{n,k}/G)\to \Brnr^0(\SL_{n,\R}/G)$, and the latter group is zero by Proposition \ref{kunyavskiitrivial}. Moreover $\Br_1^0(\SL_{n,k}/G)\simeq\rmH^1(k,\widehat{G}^\mathrm{ab})=0$ since $G^\mathrm{ab}=0$, and hence $\Brnr^0(\SL_{n,k}/G)=\Br_{\mathrm{nr},1}^0(\SL_{n,k}/G)=0$. Thus, as a corollary of Proposition \ref{kunyavskiitrivial}, we conclude that Colliot-Thélène's Conjecture \ref{conjectureColliot} predicts that the Grunwald problem for $G$ over the field $k$ admits a positive answer without any restriction on the set of places.   
\end{rmk}
\bibliographystyle{amsalpha}
\bibliography{main}
\end{document}